\documentclass[a4paper,11pt,twoside,dvips]{article}
\usepackage{latexsym}
\usepackage{amsmath}
\usepackage{amssymb}
\usepackage{textcomp} 
\usepackage[ansinew]{inputenc} 
\usepackage[dvips]{graphicx} 
\usepackage{tabularx}  
\usepackage{amsfonts}
\usepackage{float}
\usepackage{amssymb}
\usepackage{amsthm}
\usepackage{verbatim}
\usepackage{graphicx}
\usepackage{bbm}
\usepackage[T1]{fontenc}
\usepackage{fancyhdr}
\newtheorem{thm}{Theorem}[section]
\newtheorem{prop}[thm]{Proposition}
\newtheorem{lemma}[thm]{Lemma}

\newtheorem{rem}[thm]{Remark}

\newcommand{\eps}{\varepsilon}
\renewcommand{\cal}{\mathcal}
\renewcommand{\rm}{\mathrm} 
\renewcommand{\bf}{\mathbf} 
\newcommand{\bb}{\mathbb} 
\newcommand{\bbm}{\mathbbm} 

\newcommand{\df}{\stackrel{\mathrm{def}}{=}}

\numberwithin{equation}{section}

\begin{document}

\title{\textbf{Cut Points and Diffusions in Random Environment}}
\author{by Ivan del Tenno\footnote{Department of Mathematics, ETH Zurich, CH-8092 Zurich, Switzerland}\\ ETH Zurich\\
e-mail: deltenno@math.ethz.ch\\
running head: Cut Points and Diffusions in Random Environment\\[11pt]
May, 2008}
\date{\empty}
\maketitle
\begin{abstract}
In this article we investigate the asymptotic behavior of a new class of multi-dimensional diffusions in random environment. 
We introduce cut times in the spirit of the work done by Bolthausen, Sznitman and Zeitouni, see \cite{BSZ}, 
in the discrete setting providing a decoupling effect in the process. This allows us to take advantage of an ergodic structure to derive a strong law of large numbers with possibly 
vanishing limiting velocity and a central limit theorem under the quenched measure.
\end{abstract}
{\textbf{Keywords:}} cut points, diffusions in random environment, quenched invariance principle, law of large numbers, diffusive behavior

\section{Introduction}\label{a}

The object of the present article is to introduce a special class of multi-dimensional diffusions in random 
environment for which we are able to prove a law of large numbers and a functional central limit theorem governing 
the corrections to the law of large numbers, valid for a.e. environment (a so-called {\textit{quenched}} functional central 
limit theorem). The investigation of the asymptotic behavior of multi-dimensional diffusions in random environment 
is well-known for its difficulty due to the massively non-self-adjoint character of the model, and to the rarity 
of explicitly calculable examples. A special interest of the class we introduce stems from the fact that it offers 
examples of diffusions with non-vanishing random drifts where on the one-hand  no
 invariant measure for the process of the environment viewed from the particle, absolutely continuous with respect 
to the static distribution of the random environment, is known, and where on the other hand our results 
hold without certain assumptions which guarantee condition (T) or (T') of Sznitman, see \cite{SCHM1}, \cite{GOE2} 
in the continuous set-up. Thus when the limiting velocity vanishes such examples correspond to diffusive motions where 
very few results are available,  see \cite{brickup}, and \cite{BSZ}, in the discrete set-up, or \cite{SZNZEIT} 
for diffusions in random environment. And when the limiting velocity does not vanish such examples differ from 
existing results for {\textit{quenched}} functional central limit theorems, such as in the recent \cite{berzei}, 
or \cite{RS2}, in the
 discrete set-up, for such results rely on finiteness assumptions for  moments of certain regeneration times, which to the best of 
our knowledge can only be checked through some sufficient criterion for (T) or (T'). Let us mention that at present 
it is an open problem whether ballistic behavior in dimension 2 and above implies (T) or (T'). Our class contains 
examples of ballistic motion, and we do not need to check (T) or (T') (which of course does not preclude that these 
conditions may hold in these examples).  The class we introduce here is a type of continuous 
counterpart of the class considered in 
 \cite{BSZ} in the context of random walks in random environment.  The formulas we obtain for the velocity 
are reasonably explicit and might be amenable to the construction of some further examples or counterexamples, 
 in the spirit of what was done in \cite{BSZ},
  although this is not carried out here given the length of the present work. Indeed the continuous 
set-up is more delicate than the discrete set-up, and it is by no mean routine to adapt the general strategy of 
\cite{BSZ} in the context of diffusions in random environment. For an overview of results and useful techniques 
concerning this area of research we refer to 
  \cite{SZN2}, \cite{SZN}, \cite{ZEIT}.\\

Before describing our results any further, let us first introduce the model.
We consider integers $d_1\geq 5,\ d_2\geq 1$ and $d=d_1+d_2.$ The random environment is described by a probability space
$(\Omega,\cal{A},\bb{P})$ 
and we assume the existence of a group $\{\tau_x : x\in\bb{R}^{d}\}$ of $\bb{P}$-preserving transformations  on $\Omega$ that are 
jointly measurable in $x$ and $\omega.$ On $(\Omega,\cal{A},\bb{P})$ we consider an $\bb{R}^d$-valued random variable $b(\cdot)$
with vanishing first $d_1$ 
components, that is\begin{equation}\label{77}
b(\omega)= (\underbrace{0,\ldots,0}_{d_1},b^*(\omega))\in\bb{R}^d,\mbox{ for }\omega\in\Omega,
\end{equation}and we define
\begin{equation}\label{1}
b(x,\omega)\df b(\tau_x(\omega)),\qquad\mbox{ for }x\in\bb{R}^d.
\end{equation}

We assume this function to be bounded and Lipschitz continuous,
i.e. there is a constant $\kappa>0$ such that for 
all $x,y\in\bb{R}^d, \omega\in\Omega,$
\begin{equation}\label{2}|b(x,\omega)|\leq \kappa,\quad
|b(x,\omega)-b(y,\omega)|\leq \kappa|x-y|,
\end{equation}

where $|\cdot|$ denotes the Euclidean norm in $\bb{R}^d.$ We will further assume finite range dependence for the 
environment, that is for a Borel subset $F$ of $\bb{R}^d$ we define the $\sigma$-algebra
\begin{equation}\label{4}
\cal{H}_F\overset{\mbox{\scriptsize{def}}}{=}\sigma(b(x,\omega):x\in F)
\end{equation}

and assume that there is an $R>0$ such that
\begin{equation}\label{5}
\cal{H}_A \mbox{ and } \cal{H}_B \mbox{ are independent whenever }d(A,B)>R,
\end{equation}

where $d(A,B)=\inf\{|x-y|:x\in A,\ y\in B\}.$ 
We let stand $(X_t)_{t\geq 0}$ for the canonical process on $C(\bb{R}_+,\bb{R}^{d})$ and for $\omega\in\Omega, x\in\bb
{R}^d$ we denote with $P_{x,\omega}$ the unique solution to the martingale problem attached to $x$ and
\begin{equation}\label{270}
\cal{L}^{\omega}=\frac{1}{2}\Delta+ b(\cdot,\omega)\cdot\nabla,
\end{equation}i.e. the law $P_{x,\omega}$ describes the diffusion in the environment $\omega$ starting at $x$ and is 
usually called the {\textit{quenched law.}} We write $E_{x,\omega}$ for the corresponding expectation. We endow the space 
$C(\bb{R}_+,\bb{R}^{d})$ with the Borel $\sigma$-algebra $\cal{F}$ and the canonical right-continuous 
filtration $(\cal{F}_t)_{t\geq 0}.$ For the study of the asymptotic properties of $X_.$, it is convenient
to introduce the  {\it annealed law} which is the semi-direct product measure on $\Omega\times C(\bb{R}_+,\bb{R}^{d})$
defined as
\begin{equation}\label{78}
P_x\df \bb{P}\times P_{x,\omega}.
\end{equation} 

We denote with $E_x$ the corresponding expectation.
Let us mention that the laws $P_x$ typically destroy the Markovian structure but restore a useful stationarity 
to the problem.\\

Let us now explain the purpose of this work in more detail. In the first part of this article we prove a law of large 
 numbers, see Theorem \ref{89}, namely when $d_1\geq 5,$ we show that
\begin{equation}\label{201}
P_0\mbox{-a.s.,}\qquad \frac{X_t}{t}\longrightarrow v\df {E^{{P}\times K_0}
\left[\int_0^{T^1}b\left(\chi_u^{},\omega\right)du,\ T^0=0\right]},
\quad\mbox{as }t\to\infty,
\end{equation}
with a deterministic (possibly vanishing) limiting velocity $v.$ The process $\chi_u,\ u\geq 0,$ is defined 
on an enlarged probability space, see Theorem \ref{47}, on which the notion of doubly infinite bilateral 
cut times ($T^k,\ k\in\bb{Z}$) 
for the Brownian part of the diffusion is superimposed, see (\ref{19}). The definition of cut times
 involves neither the drift nor the random environment except for the parameters $\kappa$ of (\ref{2}) 
and $R$ of (\ref{5}). 
The law of $(\omega,(\chi_u^{})_{u\geq 0})$ under the measure ${P}\times K_0$ recovers 
the {\textit{annealed}} measure $P_0,$ see (\ref{17}), (1) of Theorem \ref{47} and (\ref{190}). \\
In the second part, assuming antipodal 
symmetry in the last $d_2$ components of the drift, see (\ref{100}), and when $d_1\geq 7$ (in which case $v=0$), or
when $d_1\geq 13$ without symmetry properties, we derive a functional 
central limit theorem under the \textit{quenched law,} see Theorem \ref{200}:
\begin{equation}\label{202}\begin{array}{l}
\mbox{for }\bb{P}\mbox{-a.e. }\omega,\mbox{ under the\, measure }P_{0,\omega},\mbox{ the }C(\bb{R}_+,\bb{R}^d)\mbox{\,-valued 
random}\\
\mbox{variables }B^r_.\,\df\, r^{-1/2}\,(X_{r\cdot}-vr\cdot),\ r>0,
\mbox{ where $v$ corresponds to the}\\[2pt]
\mbox{limiting velocity in (\ref{201}), converge weakly to a Brownian motion with}\\[2pt]
\mbox{deterministic covariance matrix, as $r$ tends to infinity.}
\end{array}\end{equation}
The proofs of the above results are based on the existence of so-called cut times $T^k,\ k\in\bb{Z},$ which are defined 
in a similar spirit to \cite{BSZ} and play a role comparable to the regeneration times introduced in \cite{SZNZER}. 
The assumption $d_1\geq 5$ enables to exploit the presence of these cut times and discover a decoupling effect, 
see Proposition \ref{59}. The cut times are in essence defined as follows. In the spirit of the technique applied in 
\cite{CZ} for random walks in random environments or in \cite{SHEN2} for the continuous case, we couple our 
diffusion at each integer time $n$ with an auxiliary Bernoulli variable $\Lambda_n$ such that when $\Lambda_n=1$, 
the distribution of $X_{n+1}$ given $X_n$ does not depend on the environment. We then say that a cut time occurs at the 
integer time $n,$ if the Bernoulli variable at time 
 $n-1$ takes value 1 and the future of the Brownian part of the diffusion, which corresponds to the first $d_1$ components, 
after time $n$ stays at a distance at least $2R$ from the past before time $n-1,$ see (\ref{19}) and (\ref{22}) 
for the exact definition. Due to the finite range dependence, see (\ref{5}), we then can produce decoupling 
in our process which allows an easy comparison to a process defined on a probability space with an ergodic shift in which 
we can embed an additive functional. These considerations essentially reduce
 the proof of (\ref{201}) to an application of 
Birkhoff's Ergodic Theorem. With the help of a criterion introduced by Bolthausen and Sznitman in \cite{BS}, see Lemma 4, the {\textit{quenched}} invariance principle (\ref{202}) 
follows from the {\textit{annealed}} versions, see Theorem \ref{103} and \ref{120}, by a variance calculation which involves 
a certain control on the intersections of two independent paths. 
The main strategy behind the proofs of the {\textit{annealed}} 
central limit theorems is to show an {\textit{annealed}} central limit theorem 
for a process defined as the polygonal interpolation of an ergodic process 
$Z^s_k, k\in\bb{Z},$ see (\ref{43}) and Proposition \ref{88}, which is then rescaled in time and space analogously to the 
definition of $B^n_.$ in (\ref{202}) for integers $n\geq 1,$ and which is
comparable to the original diffusion $X_.,$ see Lemma \ref{510} and (\ref{122}). The proof without symmetry 
assumption on the drift but $d_1\geq 13$ is more involved and needs an adaptation of Gordin's method, see 
for instance the proof of Theorem 7.6 in \cite{D}.\\

Let us mention that the application of Girsanov's formula yields a very handy and reasonably explicit version
of the transition density for the last $d_2$ components of the diffusion in a fixed environment given the 
Brownian part (first $d_1$ components), see (\ref{255}). The formula (\ref{255}) involves the Brownian transition 
density and the 
bridge measure which depend neither on the environment nor on the 
first $d_1$ components of the diffusion and hence enables to inspect 
the {\textit{quenched}} transition density directly. This formula is not available anymore if one wants to treat more 
general diffusions in random environment where the diffusion matrix in the last $d_2$ components becomes 
a genuinely environment dependent stationary process. Other methods would be required in this set-up, possibly in the spirit 
of filtering theory.\\

Let us now explain how this article is organized. In Section \ref{b} we couple our diffusion with a suitable 
sequence of i.i.d. Bernoulli variables, see Theorem \ref{47}. We then define the cut times $T^k,\ k\in\bb{Z},$ 
see (\ref{19}) and (\ref{22}), and provide the crucial decoupling, see Proposition \ref{59}. Finally we prove a law 
of large numbers. Section \ref{c} is dedicated to two central limit 
theorems under the annealed measure that are also consequences of the decoupling technique discussed in Section \ref{b}. 
The first central limit theorem is proved under a symmetry assumption on the drift and $d_1\geq 7,$
see (\ref{100}), whereas for the second central limit theorem $d_1\geq 13$ is assumed. In Section \ref{d} 
we show how one can strengthen the results of Section \ref{c} into central limit theorems under the quenched measure.
Finally, in the Appendix, two multidimensional versions of central limit theorems for martingales are proved. \\

{\textbf{Convention on constants:}} Unless otherwise stated, constants only depend on the quantities $d_1,d_2,\kappa,R.$ 
In calculations, generic constants are denoted by $c$ and may change from line to line, whereas $c_1, c_2, \ldots$ are constants with fixed values 
at their first appearance. With $c(q,\eta)$ we denote constants that depend on the usual parameters $d_1,d_2,\kappa,R$ 
and additionally on $q$ and $\eta.$\\

{\textbf{Acknowledgements:}} I would like to express my sincere gratitude to my advisor Prof. A.-S. Sznitman for his 
support during this work. I want also to thank T. Schmitz, L. Goergen and D. Windisch for many helpful and encouraging 
discussions.

\section{Decoupling and a law of large numbers}\label{b}

In this section we will first take advantage of the special structure of the diffusions considered in our model to couple 
auxiliary i.i.d. Bernoulli variables $\Lambda_n$ with the diffusion, see Theorem \ref{47}. Under the coupled measure, the 
distribution of the diffusion at integer time $n$ will only depend on the position at time $n-1$ and not on the 
environment when $\Lambda_{n-1}=1.$ Due to the finite range dependence, see (\ref{5}), we then discover 
with the help of cut times, which are introduced in Subsection \ref{241}, a certain decoupling effect 
under the annealed law, see Proposition \ref{59}. This finally leads to a law of large numbers, 
see Theorem \ref{89}.\\


For a real number $u\in\bb{R}$ we define its integer part as \begin{equation}\label{300}
[u]\df\sup\{n\in\bb{Z}\ |\ n\leq u\}.\end{equation}
Further we denote the $d_2$-dimensional, closed ball of radius $r>0$ centered at $y\in\bb{R}^{d_2}$ with $B^{d_2}_r(y)$ and 
write $vol(d_2)$ for its volume. For $n\geq 1, z, z'\in\bb{R}^n$ and $s>0$ we introduce the $n$-dimensional Gaussian 
kernel\begin{equation}\label{70}
p_{n}(s,z,z')\df \frac{1}{(2\pi s)^{n/2}}\exp\{-|z-z'|^2/2s\}.
\end{equation}
We denote with $W^{d_1}_0$ the set of all continuous $\bb{R}^{d_1}$-valued 
functions on $\bb{R}$ that vanish at $0.$ Furthermore we consider the space  
$W_+^{d_2}=C(\bb{R}_+,\bb{R}^{d_2})$ and the canonical coordinate processes $X^1_.,X^2_.$ defined as

\begin{equation}\begin{aligned}\label{79}
X^1_t(w)\df w(t)\mbox{ for all }t\in\bb{R}\mbox{ and }w\in W^{d_1}_0,\\
X^2_t(u)\df u(t)\mbox{ for all }t\geq 0\mbox{ and }u\in W^{d_2}_+.\hspace{0.2cm}
\end{aligned}\end{equation}

We endow the space $W^{d_1}_0$ with the $\sigma$-algebra $\cal{W}_0=\sigma(X^1_s, s\in\bb{R})$ and $W^{d_2}_+$ with the 
$\sigma$-algebra $\cal{U}=\sigma(X^2_s, s\in\bb{R}_+)$ and the canonical filtration $\cal{U}_t=\sigma(X^2_s,0\leq s\leq t),t\geq 0$
which is neither right-continuous nor complete in 
opposition to ${\cal{F}}_t,$ see above (\ref{78}).
$\bar{P}$ denotes the two-sided Wiener measure on $(W^{d_1}_0,\cal{W}_0)$ with $\bar{P}[X_0^1=0]=1.$ 
We write $\bar{E}$ for the expectation with respect to the measure $\bar{P}.$ On the measurable space 
$(W^{d_2}_+,\cal{U})$ we introduce for $y,y'\in\bb{R}^{d_2}$ the Wiener measure $\tilde{P}_y$ with $\tilde{P}_y[X^2_0=y]=1$ 
and $\tilde{P}_{y,y'},$ the Brownian bridge measure from $y$ to $y'$ on [0,1]. We write $\tilde{E}_y$ and 
$\tilde{E}_{y,y'}$ for the corresponding expectations. On the product space $(W^{d_1}_0\times W^{d_2}_+,\cal{W}_0
\otimes\cal{U})$ we define the $\bb{R}^d$-valued process
\begin{equation}\label{81}
\chi_t^{}\df (X_t^1,X_t^2),\quad t\geq 0.
\end{equation}
For $y\in\bb{R}^{d_2},\omega\in\Omega$ and $w\in W^{d_1}_0$ we denote with $\bar{K}_{y,\omega}(w)$ the probability kernel from 
$(W^{d_1}_0,\cal{W}_0)$ to $(W^{d_2}_+,\cal{U})$ defined as the unique solution of the martingale problem starting at time 0 
from $y$ and attached to
\begin{equation*}
\cal{L}_t^{w,\omega}=\frac{1}{2}\sum_{i=1}^{d_2}\partial_{ii}^2+\sum_{i=1}^{d_2}b_i^*\big((w(t),\cdot),\omega\big)
\partial_{i},
\end{equation*}see Theorem 6.3.4 in \cite{SV}. For $w\in W^{d_1}_0$ and $\omega\in\Omega$ we define the stochastic exponential 
\begin{equation}\label{13}
\cal{E}(w,\omega)\df \exp\left\{
\int_0^1 b^*((w(s),X^2_s),\omega)dX^2_s-\frac{1}{2}\int_0^1|b^*((w(s),X_s^2),\omega)|^2ds\right\},
\end{equation}which is in $L^1(\tilde{P}_{y,y'})$ for all $y,y'\in\bb{R}^{d_2},$ see the proof of 
Theorem 4.1 of \cite{LZ}, and introduce the transition density
\begin{equation}\label{255}
p_{w,\omega}(1,y,y')\df p_{d_2}(1,y,y')\tilde{E}_{y,y'}\left[\cal{E}(w,\omega)\right],
\end{equation}which is a measurable function of $\omega\in\Omega,w\in W^{d_1}_0$ and $y,y'\in\bb{R}^{d_2},$
 see for instance Theorem 44 in \cite{P} on page 158, fulfilling
\begin{equation*}
\bar{K}_{y,\omega}(w)[X^2_1\in G]=\int_{G}dy'p_{w,\omega}(1,y,y'),
\end{equation*}for all Borel sets $G$ in $\bb{R}^{d_2},$ see equation (6.35) in Chapter 5 of \cite{KS} and 
Girsanov's Formula in Theorem 6.4.2 of \cite{SV}.

\subsection{The coupling construction}
We are going to enlarge the probability space $(W^{d_1}_0\times W^{d_2}_+,\cal{W}_0\otimes\cal{U},\bar{P}\times\bar{K}_
{y,\omega})$ and provide a coupling of the process $\chi_.^{}$ with a sequence of i.i.d. Bernoulli variables, 
see Theorem \ref{47}. Let us begin with an easy fact about the
transition density defined in (\ref{255}), which will be crucial in the construction of our coupling.

{\lemma\label{91}{
Under the assumptions (\ref{77}) and (\ref{2}) on the drift $b(\cdot)$ there is a constant $\varepsilon\in (0,1)$ such that 
for all $\omega\in\Omega, w\in W^{d_1}_0, y\in\bb{R}^{d_2}$ and  $y'\in B_1^{d_2}(y)$ the following holds:
\begin{equation}\label{92}
p_{w,\omega}(1,y,y')>\frac{2\varepsilon}{vol(d_2)}.
\end{equation}}}

{\textbf{Proof:} By Jensen's inequality and (\ref{2}) we obtain that $p_{\omega,w}(1,y,y')$ is larger or equal than
\begin{equation}\label{93}
e^{-\kappa^2/2}p_{d_2}(1,y,y')\exp\left\{\tilde{E}_{y,y'}\left[\int_0^1 b^*\big((w(s),X^2_s),\omega\big)dX^2_s
\right]\right\}.
\end{equation}

Note that under the measure $\tilde{P}_{y,y'}$ the process $X_t^2, t\in [0,1],$ is a Brownian bridge from $y$ to $y'$ in time 1 and 
therefore satisfies the following stochastic differential equation, see Ex. 5.6.17 (i) and p. 354 in \cite{KS},
\begin{equation}{\begin{cases}
dX^{2}_t  =  d\beta_t+\frac{y'-X^{2}_t}{1-t}dt,&\quad0\leq t<1;\\
X^2_0=y,\ \tilde{P}_{y,y'}\mbox{-a.s. },\\
\end{cases}}
\end{equation}

for a $d_2$-dimensional standard Brownian motion $\beta_.$.
Thus, 
\begin{eqnarray*}
\left|\tilde{E}_{y,y'}\left[\int_0^1 b^*\big((w(s),X^2_s),\omega\big)dX^2_s\right]\right| & = & 
\left|\tilde{E}_{y,y'}\left[\int_0^1 b^*\big((w(s),X_s^2),\omega\big)\frac{y'-X_s^2}{1-s}ds\right]\right|\\
& \leq & \tilde{E}_{y,y'}\left[\kappa\int_0^1|\nabla_{X_s^2} \log p_{d_2}(1-s,X_s^2,y')|ds\right]\\
& \leq & c\exp\{-|y-y'|^2/c\}p_{d_2}^{-1}(1,y,y')\df g_{d_2}(y-y'),
\end{eqnarray*}

where the last inequality follows from a result of \cite{LZ}, see Theorem 2.4.
It is obvious that
\begin{eqnarray*}B^{d_2}_1(0) \ni z \longmapsto g_{d_2}(z) 
\end{eqnarray*}

is a bounded map and thus (\ref{93}) is bounded away from 0 for all $y'\in B^{d_2}_1(y).$ This finishes the proof.
\begin{flushright}$\Box$\end{flushright}

Before providing the construction of the coupling, let us introduce some further notation. We denote with
$\Lambda_.=(\Lambda_j)_{j\in\bb{Z}}$ the canonical coordinate process on $\{0,1\}^{\bb{Z}}$ and with $\cal{S}$ the canonical
product $\sigma$-algebra generated by $\Lambda_.$. We write $\lambda_.=(\lambda_j)_{j\in\bb{Z}}$ for an element of 
$\{0,1\}^{\bb{Z}}$ and with $\bf{\Lambda}^{\varepsilon},$
where $\varepsilon$ comes from (\ref{92}), we denote the unique 
probability measure on $(\{0,1\}^{\bb{Z}},\cal{S})$ under which $\Lambda_.$ becomes a sequence of i.i.d. Bernoulli 
random variables with success parameter $\varepsilon.$
We also introduce the shift operators $\{\theta_m:m\in\bb{Z}\}$ and $\{s_t:t\geq 0\}$
operating on $(W^{d_1}_0\times\{0,1\}^{\bb{Z}},\cal{W}_0\otimes\cal{S})$ and $(W^{d_2}_+,\cal{U})$ respectively such that 
\begin{eqnarray}
\label{16}\theta_m(w,\lambda_.)&=&(w(m+\cdot)-w(m),\lambda_{m+\cdot}),\\
\label{256}s_t(u)&=&u(t+\cdot).
\end{eqnarray}
Note that the pair $(w,\lambda_.)\in W^{d_1}_0\times\{0,1\}^{\bb{Z}}$ stands for the pair of processes $((w(t))_{t\in\bb{R}},$ $(\lambda_j)_{j\in\bb{Z}})$ 
whose parameter sets are different. On the product space $(W^{d_1}_0\times\{0,1\}^{\bb{Z}},\cal{W}_0\otimes\cal{S})$ 
we define the product measure 
\begin{equation}\label{17}
P\overset{\mbox{\scriptsize{def}}}{=} \bar{P}{\otimes}{\bf{\Lambda}}^{\varepsilon},
\end{equation}recall that $\bar{P}$ denotes the two-sided Wiener measure on $W^{d_1}_0$ 
 with $\bar{P}[X^1_0=0]=1.$

\pagebreak
 
{\thm\label{47}{There exists a probability kernel from $\bb{R}^{d_2}\times\Omega\times W^{d_1}_0\times\{0,1\}^{\bb{Z}}$ to 
$W^{d_2}_+$ which we denote with $K_{y,\omega}(w,\lambda_.)[O]$ for $y\in\bb{R}^{d_2},\omega\in\Omega,w\in W^{d_1}_0, 
\lambda_.\in\{0,1\}^{\bb{Z}}$ and $O\in\cal{U},$ such that:\begin{itemize}
\item[(1)] For $(y,\omega,w,\lambda_.)\in\bb{R}^{d_2}\times\Omega\times W^{d_1}_0\times\{0,1\}^{\bb{Z}},$ 
under the measure $P\times K_{y,\omega}(w,\lambda_.),$ the process $(\chi_t^{})_{t\geq 0}$ is $P_{(0,y),\omega}$-distributed, 
where $(0,y)\in\bb{R}^d,$ and in particular $W_t$ defined by
\begin{equation}\label{150}
W_t\df \chi_t^{}-(0,y)-\int_0^t b\left(\chi_s,\omega\right)ds,\mbox{ for all }t\geq 0,
\end{equation}
is a $d$-dimensional standard Brownian motion in its own filtration on 
$W^{d_1}_0\times\{0,1\}^{\bb{Z}}\times W^{d_2}_+$ endowed with the probability $P\times K_{y,\omega}(w,\lambda_.).$
\item[(2)] For each integer $n\geq 0,\ (y,\omega,w,\lambda_.)\in\bb{R}^{d_2}\times\Omega\times W^{d_1}_0\times\{0,1\}^{\bb{Z}}$ 
and any bounded measurable function $f$ on $W^{d_2}_+,$ \begin{equation}\label{257}
K_{y,\omega}(w,\lambda_.)\mbox{-a.s., }\ E^{K_{y,\omega}(w,\lambda_.)}\left[f(X^2_.)\circ s_n\,|\,\cal{U}_n\right]=
E^{K_{X^2_n,\hat{\omega}}(\theta_n(w,\lambda_.))}\left[f(X^2_.)\right],
\end{equation}with $\hat{\omega}=\tau_{(w(n),0)}(\omega).$ Moreover, \begin{equation}\label{258}
E^{K_{y,\omega}(w,\lambda_.)}\left[f(X^2_.)\right]=E^{K_{0,\tilde{\omega}}(w,\lambda_.)}\left[f(y+X^2_.)\right],
\end{equation}with $\tilde{\omega}=\tau_{(0,y)}(\omega).$
\item[(3)] For each $(y,\omega,w,\lambda_.)\in \bb{R}^{d_2}\times\Omega\times W^{d_1}_0\times\{0,1\}^{\bb{Z}}$ with 
$\lambda_0=1,$ we have that under the probability measure $K_{y,\omega}(w,\lambda_.),\ X^2_1$ is uniformly distributed on 
the ball $B^{d_2}_1(y).$ 
\item[(4)] For each integer $n\geq 0,\ (z_1,z_2)\in\bb{R}^{d_1}\times\bb{R}^{d_2},\ (y,\omega,w,\lambda_.)
\in\bb{R}^{d_2}\times\Omega\times W^{d_1}_0\times\{0,1\}^{\bb{Z}}$ and any bounded measurable function $f$ on 
$W_+^{d_2}\times C(\bb{R}_+,\bb{R}^{d}),$ we have that for $\bar{\omega}=\tau_{(z_1,z_2)}(\omega),$ 
\begin{equation}\label{75}
E^{K_{y,\bar{\omega}}(w,\lambda_.)}\left[f\bigg((X_{\cdot}^{2},b(\chi_{\cdot}^{},\bar{\omega}))_{\cdot\wedge n}\bigg)\right]
\mbox{ is }\cal{H}_{(z_1+w([0,n]))\times\bb{R}^{d_2}}\mbox{-measurable}.
\end{equation}
\end{itemize}}}In order to shorten notation we will usually not write explicitly the dependence of the kernels 
on $(w,\lambda_.)\in W^{d_1}_0\times\{0,1\}^{\bb{Z}},$ i.e. for $y\in\bb{R}^{d_2}, \omega\in\Omega$ 
we write $K_{y,\omega}$ instead of $K_{y,\omega}(w,\lambda_.).$ 
In this sense for a fixed $y\in\bb{R}^{d_2}$ we define the {\textit{annealed kernel}} from 
$(W^{d_1}_0\times\{0,1\}^{\bb{Z}},\cal{W}_0\otimes \cal{S})$ 
to $(W^{d_2}_+\times \Omega,\cal{U}\otimes{\cal{A}})$ by
\begin{equation}\label{190}K_y\df\bb{P}\times K_{y,\omega}.\end{equation}

\textbf{Proof:} Given a probability kernel $K^{(\lambda)}_{y,\omega}(w)[O]$ for $O\in\cal{U}_1, w\in W^{d_1}_0, \lambda\in
\{0,1\}, y\in\bb{R}^{d_2}$ and $\omega\in\Omega,$ which will be specified in (\ref{259}) below, there is a unique 
probability measure $K_{y,\omega}(w,\lambda_.)$ on 
$\cal{U}$ for $w\in W^{d_1}, \lambda_.\in\{0,1\}^{\bb{Z}}, y\in\bb{R}^{d_2}$ and $\omega\in\Omega$ such that for integer $m\geq 1, 
O\in\cal{U}_1:$\begin{equation}\label{650}
K_{y,\omega}(w,\lambda_.)\mbox{-a.s., }\quad K_{y,\omega}(w,\lambda_.)\left[s_m^{-1}(O)\,|\,\cal{U}_m\right]=
K_{X^2_m,\tau_{(w(m),0)}(\omega)}^{(\lambda_m)}(w(m+\cdot)-w(m))\left[O\right].
\end{equation}
An application of Girsanov's Theorem, see for instance Theorem 6.4.2 of \cite{SV}, and equation (6.35) in 
Chapter 5 of \cite{KS}, show that for 
$w\in W^{d_1}_0, y\in\bb{R}^{d_2}, \omega\in\Omega$ and $O\in\cal{U}_1,$ see below (\ref{81}) and (\ref{255}),
\begin{eqnarray*}
\bar{K}_{y,\omega}(w)\left[O\right]=\tilde{E}_y\left[\cal{E}(w,\omega),O\right]=\int_{\bb{R}^{d_2}}dy' p_{w,\omega}(1,y,y')
\frac{\tilde{E}_{y,y'}\left[\cal{E}(w,\omega), O\right]}{\tilde{E}_{y,y'}\left[\cal{E}(w,\omega)\right]}, 
\end{eqnarray*} and so, we define for $\lambda\in\{0,1\}$ and $y'\in\bb{R}^{d_2},$
\begin{equation}\label{28}{h(w,\lambda,y,y',\omega)\df\begin{cases}
\frac{\bbm{1}{\{y'\in B^{d_2}_1(y)\}}}{vol(d_2)}, 
&\mbox{if } \lambda=1;\\
\frac{1}{1-\varepsilon}\Big(p_{w,\omega}(1,y,y')-\varepsilon\frac{\bbm{1}{\{y'\in B_1^{d_2}(y)\}}}
{vol(d_2)}\Big),
&\mbox{if } \lambda=0,\\
\end{cases}}
\end{equation} and set \begin{equation}\label{259}
K_{y,\omega}^{(\lambda)}(w)\left[O\right]=\int_{\bb{R}^{d_2}}dy'h(w,\lambda,y,y',\omega)\frac{\tilde{E}_{y,y'}\left[
\cal{E}(w,\omega), O\right]}{\tilde{E}_{y,y'}\left[
\cal{E}(w,\omega)\right]}.
\end{equation}
In view of (\ref{92}), this kernel is well defined. To check the measurability of the kernel one uses a 
result of \cite{P}, see Theorem 44 on page 158. The same result can also be used to show (\ref{75}). It is then 
straightforward to see that the resulting kernel $K_{y,\omega}(w,\lambda_.)$ fulfills (1)-(4).
\begin{flushright}$\Box$\end{flushright}
{\rem{In the notation ${\bf{\Lambda}}^{\varepsilon,\lambda}[\,\cdot\,]\df{\bf{\Lambda}}^{\varepsilon}
[\,\cdot\,|\,\Lambda_0=\lambda]$ for 
$\lambda\in\{0,1\}$ and $K^{n}_{y,\omega}(w,\lambda_.)\df K_{y,\tau_{(w(n),0)}(\omega)}(w(n+\cdot)-w(n),\lambda_.)$ for 
integer $n\geq 0,\ (y,\omega,w,\lambda_.)\in\bb{R}^{d_2}\times\Omega\times W^{d_1}_0\times \{0,1\}^{\bb{Z}}$ we find as a 
consequence of (\ref{257}) and the fact that $K_{y,\omega}[X^2_{\cdot\wedge n}\in\star\,]$ depends on 
$w([0,n]),\lambda_0,\ldots,\lambda_{n-1}$ 
only, that for a fixed Brownian path $w\in W^{d_1}_0,\ {\bf{\Lambda}}^{\varepsilon}\times K_{y,\omega}\mbox{-a.s.,}$
\begin{equation}
{\bf{\Lambda}}^{\varepsilon}\times K_{y,\omega}\left[(X^2_{n+\cdot},\Lambda_{n+\cdot})\in\star\,|\,\cal{U}_n
{\otimes}\sigma(\Lambda_0,\ldots,\Lambda_{n})\right]=
{\bf{\Lambda}}^{\varepsilon,\Lambda_n}\times K^n_{X^2_n,\omega}\left[(X^2_{\cdot},\Lambda_{\cdot})\in\star\,\right]
\end{equation}}}

{\rem{\label{034}
Thank to Girsanov's Theorem we were able to construct the above kernels quite explicitly, so that we have a very 
concrete way to write expectations of $X^2_k,$ for integers $k\geq 1,$ under the quenched kernel $K_{0,\omega}$ 
using the formulas for the kernels for one time unit, see (\ref{259}). Indeed, applying (\ref{257}) successively for 
$n=k-1,\ldots,1,$ we find that for all $(w,\lambda_.)\in W^{d_1}_0\times \{0,1\}^{\bb{Z}},$\begin{equation}\label{030}
E^{K_{0,\omega}}\left[X_{k}^2\right]=E^{K_{0,\omega}}\left[E^{K_{X^2_1,\hat{\omega}_1}\circ\theta_1}\left[\cdots E^{K_{X^2_1,
\hat{\omega}_{k-1}}\circ\theta_{k-1}}\left[
X_1^2\right]\cdots\right]\right],
\end{equation}with $\hat{\omega}_i=\tau_{(w(i),0)}(\omega),\ i=1,\ldots, k-1.$ 
Using the identities (\ref{650}) and (\ref{259}) in the proof of Theorem \ref{47} 
we obtain with $y_0:=0$ that the right-hand side of (\ref{030}) equals\begin{equation}\label{031}
\int_{\bb{R}^{d_2}}\cdots\int_{\bb{R}^{d_2}}dy_1\cdots dy_{k}\prod_{i=0}^{k-1}h\left(w(i+\cdot)-w(i),
\lambda_i,y_i,y_{i+1},\hat{\omega}_i\right)y_{k}.
\end{equation} 
}}

\subsection{The cut times $T^k$}\label{241}

In this subsection we will define the cut times which are at the heart of this work, see (\ref{19}). The assumption
$d_1\geq 5$ becomes crucial at this point since it ensures the existence of these times, 
see (\ref{21}), (\ref{22}). For a concise review of the work of Erd\"os on cut times and more recent results see for instance
\cite{L2} and references therein. The consideration of these times is crucial to find certain decoupling effects in the process $\chi_.^{}$ 
under the annealed measure $P\times K_0,$ see Proposition \ref{59}, providing a comparison of $(\chi_k^{})_{k\geq 1}$ under 
$P\times K_0$ with an ergodic sequence. This enables us to deduce rather easily a law of large numbers.\\

For $r\geq 0$ and a subset $A$ of $\bb{R}^{d_1}$ we define $A^r$ as the closed $r$-neighborhood of $A.$ For 
$(w,\lambda_.)\in W^{d_1}_0\times\{0,1\}^{\bb{Z}}$ we define the set of cut times as \begin{equation}\label{19}
\cal{C}(w,\lambda_.)\overset{\mbox{\scriptsize{def}}}{=}\left\{n\in\bb{Z}\ \bigg|\ \bigg(X^1_{(-\infty,n-1]}(w)\bigg)^{R}
{\cap}\bigg(X^1_{[n,\infty)}(w)\bigg)^{R}=\emptyset,\ \Lambda_{n-1}(\lambda_.)=1\right\},
\end{equation}and consider the point process on $\bb{Z}$ \begin{equation}\label{20}
N((w,\lambda_.);dk)=\sum_{n\in\bb{Z}}\delta_n(dk)\bbm{1}_{\{n\in\cal{C}(w,\lambda_.)\}},
\end{equation}which is stationary for $\theta_1$ under the measure $P.$ It will turn out that the point process $N$ is 
double infinite, i.e. the event
\begin{equation}\label{290}
W\df\left\{(w,\lambda_.)\in W^{d_1}_0\times\{0,1\}^{\bb{Z}}\ \Big|\ N((w,\lambda_.);\bb{Z}_-)=\infty=N((w,\lambda_.);\bb{Z}_+)\right\}
\end{equation}has full $P$-probability, see Lemma \ref{87} below. We will thus restrict $P,$ see (\ref{17}), on the shift-invariant set 
$W.$ With $\cal{W}$ we denote the restriction of $\cal{W}_0{{\otimes}}\cal{S}$ to $W.$

{\rem{\label{291}On the event $\Lambda_{n-1}=1,\ n\geq 1,$ we have a very good control on the position of $\chi_n^{}$ 
by the knowledge of $\chi_{n-1}$ without any further information about the environment. Due to finite range dependence, this
will lead to a certain decoupling effect between the environment seen from the process $\chi_.^{}$ after a cut time $n$ 
and the environment affecting the process $\chi_.^{}$ before time $n-1.$ As a consequence we will find the key identity in law stated in 
Proposition \ref{59}.}}


{\lemma\label{87}($d_1\geq 5$)
\begin{eqnarray}
\label{21}&&P[0\in\cal{C}]\geq c_1(\varepsilon)>0.\\
\label{22}&&P\left[W\right]=1,\quad\mbox{and hence on }W,\quad N((w,\lambda_.);dk)=\sum_{n\in\bb{Z}}\delta_{T^n(w,\lambda_.)}(dk),
\end{eqnarray}
where $T^n,n\in\bb{Z},$ are $\bb{Z}$-valued random variables on $W$ that are increasing in $n$ such that $T^0\leq 0< T^1.$\begin{eqnarray}
\label{23}&&\hat{P}\df P[\ \cdot\ |\ 0\in\cal{C}]\mbox{ is invariant under }\hat{\theta}_1
\overset{\mbox{\scriptsize{def}}}{=}\theta_{T^1}.\\[12pt]
\label{0200}&&T^{n+m}=T^n+T^m\circ\hat{\theta}_n,\mbox{ for all }n,m\in\bb{Z}.\\[11pt]
\label{24}&&E^{\hat{P}}[T^1]=P[0\in\cal{C}]^{-1}.\\[6pt]
\label{25}&&E^{P}[f]=\frac{E^{\hat{P}}\big[\sum_{k=0}^{T^1-1}f\circ\theta_k\big]}{E^{\hat{P}}[T^1]}
\end{eqnarray}
for any bounded measurable function $f$ on $W.$
\begin{eqnarray}\label{26}
&&P[T^1>n]\leq c_2(\log{n})^{1+\frac{d_1-4}{2}}n^{-\frac{d_1-4}{2}}, \quad n\geq 1,
\end{eqnarray}
 for $c_2(\varepsilon)$ a positive constant.
}\\

\textbf{Proof:} Let us define for $w\in W^{d_1}_0,$
\begin{eqnarray*}
&&B^1_t(w)\df w(-t),\quad t\geq 0,\\
&&B^2_t(w)\df w(t),\quad t\geq 0.
\end{eqnarray*}

Noting that $B^1_., B^2_.$ and $\lambda_.$ are mutually independent and $B^1_., B_.^2$ are two $d_1$-dimensional 
standard Brownian motions on 
$(W^{d_1}_0\times\{0,1\}^{\bb{Z}},\cal{W}_0\otimes\cal{S},{P}),$ we find by using the Markov 
property of Brownian motion that
\begin{equation*}
P[0\in \cal{C}]  =  \varepsilon\int_{\bb{R}^{d_1}}p_{d_1}(1,0,x){P}\left[\left(x+B^1_{[0,\infty)}\right)^{R}\cap
\left(B^2_{[0,\infty)}\right)^{R}=\emptyset\right]dx.
\end{equation*}

To prove (\ref{21}) it suffices to show that for some set $A\subseteq\bb{R}^{d_1}$ of positive Lebesgue measure, 
 
\begin{equation}\label{31}
{P}\left[\left(x+B^1_{[0,\infty)}\right)^{R}\cap\left(B^2_{[0,\infty)}\right)^{R}=\emptyset\right]>0
\mbox{\quad for all }x\in A.
\end{equation}
For $i,j\geq 0$ let us define the event
\begin{equation}\label{32}A_{i,j}=\left\{\left(B^1_{[i,i+1]}\right)^{R}\cap\left(B^2_{[j,j+1]}\right)^{R}\ne\emptyset\right\}.\end{equation}
From the Markov property and the independence of  $B^1_.$ and $ B_.^2,$ it follows for $(i,j)\ne(0,0)$ that
\begin{equation}\begin{split}
{P}\left[A_{i,j}\right]\   = \ &\int_{\bb{R}^{d_1}}p_{d_1}(i+j,0,x)
{P}\Big[\Big(x+B^1_{[0,1]}\Big)^R\cap \Big(B^2_{[0,1]}\Big)^R\ne\emptyset\Big]dx\\
\ \leq\ & \int_{\bb{R}^{d_1}}p_{d_1}(i+j,0,x){P}
\Big[|x|\leq \sup_{0\leq s\leq 1}|B_s^1|+ \sup_{0\leq s\leq 1}|B_s^2|+2R\Big]dx.
\end{split}\end{equation} Using Fubini and the fact that $p_{d_1}(i+j,0,z)\leq c (i+j)^{-d_1/2}$ we obtain that
\begin{equation}\label{280}
P[A_{i,j}]\leq \frac{c}{(i+j)^{d_1/2}}\left(E^{P}\left[\sup_{0\leq s\leq 1}|B^1_s|^{d_1}\right]+R^{d_1}\right) \leq \frac{c}{(i+j)^{d_1/2}},
\end{equation}which implies, since $d_1\geq 5,$
\begin{equation}\label{34}
\sum_{i,j=0}^{\infty}P[A_{i,j}]<\infty.
\end{equation}


In analogy to the proof of Proposition 3.2.2 in \cite{L}, where intersection probabilities of two independent random 
walks are investigated, we call $(i,j)$ a *-last intersection if $A_{i,j}$ occurs while 
$A_{i',j'}$ for $i'\geq i, j'\geq j$ with $(i',j')\ne(i,j)$ do not. Because of (\ref{34}) and Borel-Cantelli's Lemma we know that ${P}$-a.e. 
pair of paths $(B^1_t(w))_{t\geq 0}, (B^2_t(w))_{t\geq 0}$ has at least one such *-last intersection. Hence
\begin{equation*}
1\leq\sum_{i=0}^{\infty}\sum_{j=0}^{\infty}{P}\left[(i,j)\mbox{ is a *-last intersection}\right],
\end{equation*}
which implies the existence of a pair $(I,J)$ such that
\begin{eqnarray*}
0 & < & {P}\left[(I,J)\mbox{ is a *-last intersection}\right]\leq {P}\left[\left(B^1_{[I+1,\infty)}\right)^{R}\cap
\left(B^2_{[J+1,\infty)}\right)^{R}=\emptyset\right]\\
&= & \int_{\bb{R}^{d_1}}p_{d_1}(I+J+2,0,x){P}
\left[\left(x+B^1_{[0,\infty)}\right)^{R}\cap\left(B^2_{[0,\infty)}\right)^{R}=\emptyset\right]dx,
\end{eqnarray*}where in the last equality we used the Markov property and the independence of $B^1_.$ and $B_.^2.$ Since the integrand is non-negative, this proves 
(\ref{31}) and hence (\ref{21}). By an analogous result for simple stationary point processes on $\bb{Z}$ as 
Lemma II.12 in \cite{N}, one finds using the ergodicity of $\theta_1$ that (\ref{22})  holds true. The 
measure $\hat{P}$ corresponds up to a multiplicative constant to the Palm measure attached to the 
stationary point process $N,$ see 
Chapter II in \cite{N}, in particular (10) on page 317. The statements (\ref{23})-(\ref{25}) 
are then standard consequences. Note that (\ref{25}) is a consequence of (19) on page 331 of 
\cite{N} and that (\ref{24}) follows from (\ref{25}) with the choice $f=\bbm{1}_{\{0\in\cal{C}\}}.$ 
It remains to show (\ref{26}). For integer $L\geq 1$ and for $j\geq 0,$ we define $$k_j:=1+Lj.$$ For  $J\geq 1$ we find 
that\\

$P\left[T^1>k_{3J}\right]=P\left[N\left((w,\lambda_.);[1,k_{3J}]\right)=0\right]$
\begin{eqnarray*}
&&\leq P\left[N\left((w,\lambda_.);[1,k_{3J}]\right)=0,\,\bigcap_{j=0}^{3J}\left\{\left(X^1_{(-\infty,k_j-1]}\right)^{R}\cap
\left(X^1_{[k_{j+1},\infty)}\right)^{R}=\emptyset\right\}\right]\\
&&+\ \sum_{j=0}^{3J}P\left[\left(X^1_{(-\infty,k_j-1]}\right)^{R}\cap
\left(X^1_{[k_{j+1},\infty)}\right)^{R}\ne\emptyset\right]\\[11pt]
&&=:a_1+a_2.
\end{eqnarray*}
First we bound $a_2.$ Note that for integer $n\geq 1,$\begin{equation}\label{33}
P\left[\left(B^1_{[0,\infty)}\right)^{R}\cap\left(B^2_{[n,\infty)}\right)^{R}\ne\emptyset\right]\leq 
\sum_{i\geq 0,j\geq n}P[A_{i,j}] \overset{(\ref{280})}{\leq} c\sum_{j\geq n}j^{1-\frac{d_1}{2}}\leq
c n^{-\frac{d_1-4}{2}},
\end{equation}and hence by stationarity of Brownian motion,\begin{equation}\label{087}
a_2=(3J+1)P\left[\left(B^1_{[0,\infty)}\right)^{R}\cap\left(B^2_{[L+1,\infty)}\right)^{R}\ne\emptyset\right]\leq
c(3J+1)(L+1)^{-\frac{d_1-4}{2}}.
\end{equation}
Now we turn to the control of $a_1.$ For $j=1,\ldots,3J,$ observe that on the event $\{N((w,\lambda_.);[1,k_{3J}])=0)\},$ the 
following inclusion holds:
\begin{eqnarray*}
&\left\{\left(X^1_{(-\infty,k_{j-1}-1]}\right)^{R}\cap\left(X^1_{[k_{j},\infty)}\right)^{R}=\emptyset\right\}
\cap
\left\{\left(X^1_{(-\infty,k_j-1]}\right)^{R}\cap\left(X^1_{[k_{j+1},\infty)}\right)^{R}=\emptyset\right\}&\\[10pt]
&{\displaystyle{\subseteq}} \left\{\left(X^1_{[k_{j-1}-1,k_j-1]}\right)^{R}\cap\left(X^1_{[k_{j},k_{j+1}]}\right)^{R}\ne\emptyset\right\}
\cup\bigg\{\lambda_{k_j-1}=0\bigg\}.&
\end{eqnarray*}
We thus find that the event 
\begin{equation*}
\bigcap_{j=3,6,\ldots}^{3J}\left\{\left(X^1_{[k_{j-1}-1,k_j-1]}\right)^{R}\cap\left(X^1_{[k_{j},k_{j+1}]}\right)^{R}
\ne\emptyset\right\}
\cup\bigg\{\lambda_{k_j-1}=0\bigg\}
\end{equation*}
occurs, whenever the event considered in $a$ occurs. By independence of Brownian increments and the fact that 
$\theta_1$ preserves $P$ we obtain that
\begin{equation}\label{35}
a_1\leq P\left[\left\{\left(X^1_{[0,L]}\right)^{R}\cap\left(X^1_{[L+1,2L+1]}\right)^{R}\ne\emptyset\right\}
\cup\bigg\{\lambda_{L}=0\bigg\}\right]^J\leq P\left[0\notin\cal{C}\right]^J.
\end{equation}
Choosing a large enough $\gamma$ which depends on $d_1, R$ and $ \varepsilon$ and setting $J=[\gamma \log n],
L=[\frac{n}{3J}],$ we obtain (\ref{26}) from (\ref{087}) and (\ref{35}).\begin{flushright}$\Box$\end{flushright}

\subsection{A decoupling effect and a law of large numbers}\label{decoupling}
         
Now we will exploit the presence of cut times, see (\ref{22}), in order to produce  decoupling in the 
process $\chi_.^{}$ under the measure $P\times K_0,$ see (\ref{190}). For this purpose we introduce the process $Z_.$ living on an 
enlarged space, see below (\ref{36}), equipped with a measure $Q^0,$ that uses our previous 
coupling construction and the cut times, see (\ref{37}) and (\ref{38}). The idea behind the construction of the 
process $Z_.$ is to start after each cut time a fresh path for $X^2_.$ in a new environment, which is chosen 
independently from the previous environment, see Remark \ref{291}. We then recover the law of the process $\chi_.^{}$ 
at integer times under $P\times K_0,$ see Proposition \ref{59}.\\

First we have to introduce some further notation. Consider the product spaces
\begin{equation}\label{36}
\Gamma^0\df W\times(W_+^{d_2}\times\Omega)^{\bb{N}},\qquad
\Gamma^s\df W\times(W_+^{d_2}\times\Omega)^{\bb{Z}}
\end{equation}
endowed with their product $\sigma$-algebras, see (\ref{290}) for the definition of $W.$ Recall at this point the definition of $\hat{P},$ see 
(\ref{23}), and note that in the sequel all the measures denoted with a $\hat{\ }$ correspond up to a different normalization to the Palm measure attached to the point
process $N((w,\lambda_.);dk),$ see (\ref{20}). On the spaces defined in (\ref{36}) we introduce the measures
\begin{equation}\label{37}
Q^0\df P\times M^0,\qquad \hat{Q}^0\df \hat{P}\times M^0 ,\qquad Q^s\df P\times M^s,\qquad \hat{Q}^s\df \hat{P}\times M^s,
\end{equation}
where $M^0$ and $M^s$ stand for the kernels from $W$ to $(W^{d_2}_{+}\times\Omega)^{\bb{N}}$ respectively from
$W$ to $(W^{d_2}_{+}\times\Omega)^{\bb{Z}}$ defined by
\begin{equation}\label{38}
M^0((w,\lambda_.);d\gamma^0)=K_0((w,\lambda_.);du_0d\omega_0){{\otimes}}\bigotimes_{m\geq 1}
K_0(\theta_{T^m}(w,\lambda_.);du_md\omega_m),
\end{equation}recall the definition (\ref{190}),
with $\gamma^0=(u_m,\omega_m)_{m\geq 0}\in (W_+^{d_2}\times\Omega)^{\bb{N}},$ and similarly
\begin{equation}\label{39}
M^s((w,\lambda_.);d\gamma^s)=\bigotimes_{m\in\bb{Z}}
K_0(\theta_{T^m}(w,\lambda_.);du_md\omega_m)
\end{equation}
with $\gamma^s=(u_m,\omega_m)_{m\in\bb{Z}}\in (W_+^{d_2}\times\Omega)^{\bb{Z}}.$ On $\Gamma^0$ we define the process 
$(Z_t)_{t\geq 0}$ by
\begin{equation}\label{40}
Z_t\df (X_t^1,Y_t), t\geq 0,
\end{equation}
with $X^1_.$ defined in (\ref{79}) and   
\begin{equation}\label{41}
\begin{aligned}
Y_t & \df & u_0(t),\hspace{3.3cm}\quad\mbox{ for }0\leq t<T^1,\mbox{ and }\\
Y_{(T^m+t)\wedge T^{m+1}} & \df & Y_{T^m}+u_{m}(t\wedge (T^{m+1}-T^m)),\quad\mbox{ for }m\geq 1,t\geq 0.
\end{aligned}
\end{equation}Note that $Z_0=0,\ Q^0$-a.s.. Loosely speaking, the process $Z_.$ is constructed by attaching after each cut time
a new path for the $\bb{R}^{d_2}$-components which evolves in a new 
independent environment. 
Similarly we define the two-sided process $(Z^s_t)_{t\in\bb{R}}$ on $\Gamma^s$ by
\begin{equation}\label{43}
Z_t^s\df (X^1_t,Y_t^s), t\in\bb{R},
\end{equation}
where for $m\in\bb{Z}, t\in\bb{R}_+,$
\begin{equation}\label{44}\begin{array}{rcl}
Y^s_0&\df &0,\\
Y^s_{(T^m+t)\wedge T^{m+1}}&\df& Y_{T^m}^s+u_m(t\wedge (T^{m+1}-T^m)),
\end{array}\end{equation}
and we introduce also the $\Omega$-valued process $(\alpha_t^s)_{t\in\bb{R}}$ by 
\begin{equation}\label{45}
\alpha_t^s\df\tau_{{\scriptscriptstyle{Z^s_t-Z^s_{T^m}}}}(\omega_m),\mbox{ for }T^m\leq t<T^{m+1},m\in\bb{Z},
\end{equation}which plays the role of the "relevant environment viewed from the particle". 
Note that by definition, $Z^s_0=0,\ Q^s$-a.s..
{\rem{\label{540}Note that by definition we have that under the measure $\hat{Q}^s,$ the joint distribution of $T^1$ 
and the piece of trajectory $(Z^s_t)_{t\in  [0,T^1]}$ 
is the same as the joint distribution of $T^1$ and $(\chi_t^{})_{t\in [0,T^1]}$ under $\hat{P}\times K_{0},$ see (\ref{79}),
(\ref{81}) for the definition of $\chi_.^{}$ and recall that $\hat{P}[T^0=0]=1.$}}\\

The following proposition yields a crucial identity in law.

{\prop\label{59}
Under the measure $Q^0$, the sequence of random vectors $(Z_n)_{n\geq 0}$ has the same law as $(\chi_n^{})_{n\geq 0}$ 
under the measure $P\times K_0.$
}\\

{\textbf{Proof:}} The idea of the proof is to fix $(w,\lambda_.)\in W$ and then to show by induction that for all integers 
$m\geq 0$ the following statement holds:
\begin{equation}\label{350}\begin{array}{c}
\mbox{For all bounded measurable functions }f^k,k=0,\ldots,m,\mbox{ on }\bb{R}^d,\\[11pt]
E^{K_0}\left[\prod_{k=0}^{m}f^k\left(\chi_k^{}\right)\right]=
E^{M^0}\left[\prod_{k=0}^{m}f^k\left(Z_k\right)\right].
\end{array}\end{equation} Proposition (\ref{59}) then follows by integrating out with respect to $P,$ see (\ref{37}) for the definition of $Q^0.$
Let us fix $(w,\lambda_.)\in W$ and note that (\ref{350}) holds true for $0\leq m\leq T^1(w,\lambda_.)$ by definition, 
see (\ref{38}),(\ref{40}) and (\ref{41}). We assume the above statement to be true for $m$ and show that it must still hold 
for $m+1.$ Without loss of generality we can assume that $l=T^N<m+1\leq T^{N+1}$ for an integer $1\leq N\leq m.$ Recall 
that $K_0=\bb{P}\times K_{0,\omega},$ see (\ref{190}), and so, applying (\ref{257}) with $n=l$ and then with $n=l-1$ we obtain 
that\begin{equation}\label{700}\begin{array}{lll}
E^{K_0}\left[\prod_{k=0}^{m+1}f^k(\chi_k^{})\right]
&=&\bb{E}\times E^{K_{0,\omega}}\bigg[
\prod_{k=0}^{l-1}f^k(\chi_k^{})E^{K_{X^2_{l-1},\hat{\omega}}\circ\theta_{l-1}}\bigg[f^l(w(l),X^2_1)\times\\[11pt]
&&E^{K_{X^2_{1},\tilde{\omega}}\circ\theta_{l}}\bigg[\prod_{k=1}^{m+1-l}f^{l+k}(w(l+k),X_k^2)\bigg]\bigg]\bigg]
\end{array}\end{equation}
with $\hat{\omega}=\tau_{(w(l-1),0)}(\omega)$ and $\tilde{\omega}=\tau_{(w(l),0)}(\omega).$ Since $l=T^N$ is a cut times, we have that $\lambda_{l-1}=1,$ see (\ref{19}), and hence with (3) of Theorem \ref{47} and (\ref{258}) we  
find that the right-hand side of (\ref{700}) is equal to\begin{equation}\label{701}\begin{array}{l}
{\displaystyle{\int_{\bb{R}^{d_2}}}}\frac{dy}{vol(d_2)}\bb{E}\bigg[
E^{K_{0,\omega}}\bigg[\prod_{k=0}^{l-1}f^k(\chi_k^{})\bbm{1}_{\{y\in B^{d_2}_1(X^2_{l-1})\}}\bigg]f^l(w(l),y)\times\\[11pt]
E^{K_{0,\bar{\omega}}\circ\theta_l}\bigg[\prod_{k=1}^{m+1-l}f^{l+k}(w(l+k),y+X^2_k)\bigg]
\bigg]
\end{array}\end{equation}with $\bar{\omega}=\tau_{(w(l),y)}(\omega),$ where we also used Fubini's Theorem. 
From the definition of the cut times $T^k$, 
see (\ref{19}) and Lemma \ref{87}, and the measurability property 
(\ref{75}), we see that all the factors under the $\bb{P}$-expectation in (\ref{701}) are independent, see (\ref{5}). 
Together with the induction hypothesis and stationarity of the environment (i.e. $\tau_x\bb{P}=\bb{P}$), 
we obtain that (\ref{701}) is equal to\begin{equation}\label{702}\begin{array}{l}
{\displaystyle{\int_{\bb{R}^{d_2}}}}\frac{dy}{vol(d_2)}
E^{M^0}\bigg[\prod_{k=0}^{l-1}f^k(Z_k)\bbm{1}_{\{y\in B^{d_2}_1(Y_{l-1})\}}\bigg]f^l(w(l),y)\times\\[11pt]
E^{K_{0}\circ\theta_l}\bigg[\prod_{k=1}^{m+1-l}f^{l+k}(w(l+k),y+X^2_k)\bigg].
\end{array}\end{equation}
Recalling the definitions (\ref{38}), (\ref{40}) and (\ref{41}), we deduce with the help of Fubini's 
Theorem that (\ref{702}) equals$$
E^{M^0}\bigg[\prod_{k=0}^{l}f^k(Z_k)E^{K_{0}\circ\theta_l}\bigg[\prod_{k=1}^{m+1-l}f^{l+k}(w(l+k),Y_l+X^2_k)\bigg]\bigg]
=E^{M^0}\bigg[\prod_{k=0}^{m+1}f^k(Z_k)\bigg],
$$where we used that $\theta_l=\theta_{T^N}.$ This finishes the induction step.\begin{flushright}$\Box$\end{flushright}

{\rem{By construction of the probability kernel $K_{y,\omega}(w,\lambda_.),\ 
y\in\bb{R}^d,\ \omega\in\Omega,\ (w,\lambda_.)\in W,$
see in particular (3) of Theorem \ref{47}, we have that due to $\lambda_{T^k-1}=1,k\geq 1,$ see the 
definition of cut times (\ref{19}), the transition from $X^2_{T^k-1}$ to $X^2_{T^k}$ depends only on the 
position  $X^2_{T^k-1}$ without any additional information on the environment. However, the piece of trajectory 
$X^2_t,\ T^k-1\leq t\leq T^k,$  is influenced by the environment, see (\ref{75}).
That is the reason why a decoupling effect concerning the environment, as described by the process $Z_t,t\geq 0,$
under $Q^0,$ can only be observed in the original process 
$\chi^{ }_t,t\geq 0,$ under $P\times K_0$ when we ignore the piece of trajectory during one unit of time 
just before each cut time.
In fact it can be shown by the same arguments as in the proof of Lemma \ref{59} 
that under $P\times K_0,$ the sequence of random variables 
$\chi^{ }_{(T^m+\cdot)\wedge(T^{m+1}-1)}, m\geq 0,$ has the same law as 
$Z_{(T^m+\cdot)\wedge(T^{m+1}-1)}, m\geq 0,$ under $Q^0,$ where we set $T^0=0.$ 
}}\\

We now introduce on $\Gamma^s$ a shift $(\Theta_k)_{k\in\bb{Z}}$ via:
\begin{equation}\label{51}
\Theta_k((w,\lambda_.),\gamma^s)=(\theta_k(w,\lambda_.),(u_{m+n},\omega_{m+n})_{m\in\bb{Z}})\mbox{ on }
T^n(w,\lambda_.)\leq k< T^{n+1}(w,\lambda_.),
\end{equation}
with $\gamma^s=(u_m,\omega_m)_{m\in\bb{Z}}\in (W^{d_2}_+\times\Omega)^{\bb{Z}}.$ 
{\prop{\label{88}For all $\bar{\gamma^s}\in\Gamma^s$ the following identities hold:
\begin{eqnarray}
\label{141} && Z^s_{l+a}(\bar{\gamma^s})-Z^s_l(\bar{\gamma^s})=Z^s_a\circ\Theta_l(\bar{\gamma^s}),\mbox{ for }l\in\bb{Z},
\ a\in\bb{R}_+,\\
\label{52}&& Z^s_n(\bar{\gamma^s})=\sum_{k=0}^{n-1}Z^s_1\circ\Theta_k(\bar{\gamma^s}),\mbox{ for integers }n\geq 1,\\
\label{53}&& \alpha_u^s(\bar{\gamma^s})=\alpha_{u_r}^s\circ\Theta_{[u]}(\bar{\gamma^s}),\mbox{ for }u=[u]+u_r\in\bb{R}.
\end{eqnarray}
Moreover,\begin{eqnarray}\label{54} && \Theta_1\mbox{ preserves }Q^s\mbox{ and in fact }(\Gamma^s,\Theta_1,Q^s)\mbox{ is ergodic,}\\
\label{234}&& E^{Q^s}\left[f\right]=\frac{E^{\hat{Q}^s}\big[\sum_{k=0}^{T^1-1}f\circ\Theta_k\big]}{E^{\hat{P}}\left[T^1\right]},
\mbox{ for any bounded}\\
\nonumber&&\mbox{measurable function $f$ on $\Gamma^s,$}\\[6pt]
\label{252} && Z_1^s\in L^m(Q^s)\mbox{ for all }m\in[1,\infty),\mbox{ when }d_1\geq 5.
\end{eqnarray}
}}

{\textbf{Proof:}} The identities (\ref{141})-(\ref{53}) follow by direct inspection of the definitions 
(\ref{43})-(\ref{45}) and (\ref{51}). The proof of (\ref{54}) exactly follows the 
proof in the discrete setting (see \cite{BSZ}, page 534-535). There are slight differences in the notation. 
The objects $\Gamma_s, Q_s, M_s, \{\tilde{w}_m\}_{m\in\bb{Z}}$ and $\cal{D}$ in \cite{BSZ} correspond to 
our $\Gamma^s, Q^s, M^s, \{u_m\}_{m\in\bb{Z}}$ and $\cal{C}.$ Further one has to read $(w,\lambda_.)$ 
instead of $w$ in the proof in \cite{BSZ}. Let us point out that the main strategy in showing the ergodicity 
of $(\Gamma^s,\Theta_1,Q^s)$ is to prove that $(\Gamma^s\,{{\cap}}\,\{0\in\cal{C}\},\hat{\Theta}_1
\overset{\mathrm{\scriptscriptstyle{def}}}{=}\Theta_{T^1},\hat{Q}^s)$ is ergodic, which is indeed an equivalent statement,
see (34) on page 357 in \cite{N}.
Analogously to (\ref{25}) we find (\ref{234}) as a standard consequence of the first part of the statement in (\ref{54}). 
We now come to the proof of 
(\ref{252}). We choose $m\in [1,\infty),$ 
then by definition (\ref{43}),\begin{equation}\label{310}
E^{Q^s}\left[|Z_1^s|^m\right]\leq 2^{m-1}\left\{E^{Q^s}\left[|X_1^1|^m\right]+E^{Q^s}\left[|Y_1^s|^m\right]\right\}.
\end{equation}The first expectation on the right-hand side of (\ref{310}) is finite since $X^1_.$ is a standard 
$d_1$-dimensional Brownian motion under $Q^s.$ In the notation (\ref{37}), (\ref{39}) and (\ref{44}) we have that
\begin{eqnarray*}
E^{Q^s}\left[|Y_1^s|^m\right]& = & E^{P}\left[E^{K_0\circ\theta_{T^0}}\left[|u_0(1-T^0)-u_0(-T^0)|^m\right]\right]\\
& = & \sum_{n\geq 0}  E^{P}\left[T^0=-n,\ E^{K_0\circ\theta_{-n}}\left[|u_0(1+n)-u_0(n)|^m\right]\right]\\
&\overset{stat.}{=}& \sum_{n\geq 0}  E^{\hat{P}}\left[T^1>n,\ E^{K_0}\left[|u_0(1+n)-u_0(n)|^m\right]\right]
P[T^0=0]
\end{eqnarray*}
If we show that the above expectation with respect to the measure $K_0$ is uniformly bounded, then (\ref{252}) follows since 
$\sum_{n\geq 0}\hat{P}[T^1>n]=E^{\hat{P}}[T^1]=P[T^0=0]^{-1}<\infty,$ see (\ref{24}) and (\ref{21}).
Indeed, by construction of the kernel $K_0=\bb{P}\times K_{0,\omega},$ see (\ref{190})-(\ref{259}), we find that
for each fixed $(w,\lambda_.)\in W,$\begin{eqnarray}
\nonumber && E^{K_0}\left[\left|u_0(1+n)-u_0(n)\right|^{m}\right]\\
\label{002}&&=
\bb{E}\left[E^{K_{0,\omega}}\left[\int_{\bb{R}^{d_2}}
h(w(n+\cdot)-w(n),\lambda_n,u_0(n),y,\hat{\omega})\left|y-u_0(n)\right|^{m}dy
\right]
\right],
\end{eqnarray}with $\hat{\omega}=\tau_{(w(n),0)}(\omega).$ When $\lambda_n=1,$ we immediately see by the definition of $h,$ see (\ref{28}), that the integral under 
the expectation is $K_{0,\omega}$-a.s. bounded by 1. In the other case, when $\lambda_n=0,$ the above integral is 
$K_{0,\omega}$-a.s. less or equal to\begin{equation}\label{003}
\frac{1}{1-\eps}\int_{\bb{R}^{d_2}}p_{w(n+\cdot)-w(n),\hat{\omega}}(1,u_0(n),y)\left|y-u_0(n)\right|^{m}dy+\frac{\eps}{1-\eps}.
\end{equation}A result in \cite{oleinik} concerning exponential bounds 
on fundamental solutions of parabolic equations of second order, see Theorem 1 on page 67, tells us that\begin{equation}\label{004}
p_{w(n+\cdot)-w(n),\hat{\omega}}(1,u_0(n),y)\leq c_3(w,\omega)e^{-c_4(w,\omega)|y-u_0(n)|^2},
\end{equation}for some positive constants $c_3(w,\omega),\ c_4(w,\omega).$ A closer look into the proof of the applied result from 
\cite{oleinik} reveals that the constants $c_3$ and $c_4$ in (\ref{004}) can indeed be chosen to be independent of the
Brownian path $w$ and the environment $\omega$ due to the uniform boundedness and Lipschitz constant of the the 
drift $b,$ see (\ref{2}). With  this in mind, combining (\ref{004}) and (\ref{003}) one easily sees that (\ref{002}) is also uniformly bounded in the case when $\lambda_n=0.$ This finishes the proof of (\ref{252}).\begin{flushright}$\Box$\end{flushright}

Now we are ready to state a law of large numbers when $d_1\geq 5$. For the notation see (\ref{190}), (\ref{23}), (\ref{37}), 
(\ref{39}), (\ref{43})-(\ref{45}).

{\thm{\label{89}($d_1\geq 5$)\begin{equation}\label{55}
P_0\mbox{-a.s.,}\qquad \frac{X_t}{t}\underset {t\to\infty}{\longrightarrow} v\df 
\frac{E^{\hat{P}\times K_0}
\left[\int_0^{T^1}b\left(\chi_u^{},\omega\right)du\right]}{E^{\hat{P}}[T^1]}
=E^{Q^s}\left[\int_0^1 b(\alpha_u^s)du\right]=
E^{Q^s}\left[Z_1^s\right].
\end{equation}
}}

{\textbf{Proof:}}
 First we prove that\begin{equation}\label{293} P_0\mbox{-a.s.,}\quad \lim_{t\to\infty}\frac{X_t}{t}=
E^{Q^s}\left[Z_1^s\right].\end{equation} 
For all $t\geq 1,$
\begin{equation}\label{56}
\left|\frac{X_t}{t}-E^{Q^s}\left[Z_1^s\right]\right|\leq \frac{1}{t}\left|X_t-X_{[t]}\right|
+\left|\frac{X_{[t]}}{[t]}\cdot\frac{[t]}{t}-E^{Q^s}\left[Z_1^s\right]\right|.
\end{equation}
For $\omega\in\Omega,$ under $P_{0,\omega}$ the process  $(W^{'}_t)_{t\geq 0}$ defined as 
$W'_t\df X_t-X_0-\int_0^{t}b(X_s,\omega)ds$ is a $d$-dimensional Brownian motion 
and $P_{0,\omega}$-a.s., 
\begin{equation}\label{57}
\begin{aligned}
\frac{1}{t}\left|X_t-X_{[t]}\right| & = \frac{1}{t}\left|\int_{[t]}^t b(X_s,\omega)ds
+\int_{[t]}^t dW^{'}_s\right|\\
& \leq \frac{1}{t}\left(\kappa + \left|W^{'}_t-W^{'}_{[t]}\right|\right).
\end{aligned}
\end{equation}A standard application of Borel-Cantelli's Lemma and Bernstein's inequality shows that the last expression in 
(\ref{57}) converges $P_{0,\omega}$-a.s. to 0, as $t\to\infty.$ Together with (\ref{56}) we see that to prove (\ref{293}) 
it suffices to show for integers $n\geq 1$ that $P_0$-a.s., $\frac{1}{n}X_n$ converges to 
$E^{Q^s}\left[Z_1^s\right],$ as $n\to\infty.$ As a consequence of (1) of Theorem \ref{47} and 
Proposition \ref{59}, we therefore obtain (\ref{293}), once we show that $\frac{Z_n}{n}\longrightarrow 
E^{Q^s}[Z^s_1],\ Q^0$-a.s., as $n\to\infty.$ As we will now see, the latter claim follows from the convergence of 
$\frac{Z^s_n}{n}$ under $Q^s,$ which is an immediate consequence of (\ref{52}), (\ref{54}), (\ref{252}) 
and Birkhoff's Ergodic 
Theorem. Indeed, we construct an enlarged probability space on which both processes $Z_.$ and $Z_.^s$ can be defined. 
Consider the product space
\begin{equation}\label{60}
\Gamma\df W\times\left(W^{d_2}_+\times\Omega\right)\times
\left(W^{d_2}_+\times\Omega\right)^{\bb{Z}}
\end{equation}
endowed with its product $\sigma$-algebra and the measure 
\begin{equation}\label{61}
Q\df P\times M,
\end{equation}
where $M$ is  the probability kernel from $W$ to $\left(W^{d_2}_+\times\Omega\right)\times
\left(W^{d_2}_+\times\Omega\right)^{\bb{Z}}$ defined as
\begin{equation}\label{62}
M((w,\lambda_.);d\gamma)=K_0((w,\lambda_.);du_0^{'}d\omega_0^{'}){{\otimes}}\bigotimes_{m\in\bb{Z}}
K_0(\theta_{T^m}(w,\lambda_.);du_md\omega_m)
\end{equation}
with $\gamma=((u_0^{'},\omega_0^{'}),(u_m,\omega_m)_{m\in\bb{Z}})\in (W^{d_2}_+\times\Omega)\times
(W^{d_2}_+\times\Omega)^{\bb{Z}}.$ With the projections
\begin{equation}\label{505}\begin{array}{l}
\pi^0:((w,\lambda_.),\gamma)\in\Gamma\longmapsto ((w,\lambda_.),(u_0^{'},\omega_0^{'}),(u_m,\omega_m)_{m\geq 1})
\in\Gamma^0,\\
\pi^s:((w,\lambda_.),\gamma)\in\Gamma\longmapsto ((w,\lambda_.),(u_m,\omega_m)_{m\in\bb{Z}})\in\Gamma^s,
\end{array}\end{equation}we find that $Q^0=\pi^0\circ Q$ and $Q^s=\pi^s\circ Q.$ We thus obtain that under $Q,$ the processes
\begin{equation}\label{107}
\tilde{Z}_t\df Z_t\circ\pi^0, t\geq 0,\mbox{ and }\tilde{Z}_t^s\df Z_t^s\circ\pi^s, t\in\bb{R},
\end{equation}
defined on $\Gamma$ have the same law as our original processes $Z_.$ and $Z_.^s$ under $Q^0$ and 
$Q^s$ respectively. Since $Q$-a.s.,
\begin{equation}\label{63}
\tilde{Z}_{T^1+t}-\tilde{Z}_{T^1}=\tilde{Z}^s_{T^1+t}-\tilde{Z}^s_{T^1}\mbox{\quad for all }t\geq 0,
\end{equation}
it follows that $Q$-a.s.,
\begin{equation}\label{64}
\frac{1}{t}\left|\tilde{Z}_t-\tilde{Z}_t^s\right|\leq \frac{1}{t}\sup_{a\in[0,T^1]}\left|
\tilde{Z}_a-\tilde{Z}_a^s\right|\overset{\scriptscriptstyle{t\to\infty}}{\longrightarrow} 0.
\end{equation}
We thus find that $\frac{\tilde{Z^s_n}}{n}$ and $\frac{\tilde{Z_n}}{n}$ have the same limit $Q$-a.s., which concludes the proof of 
(\ref{293}). We now show the second and the third 
equality in (\ref{55}). First we show that\begin{equation}\begin{aligned}\label{65}
\lim_{n\to\infty}\frac{E^{\hat{P}\times K_0}\left[\int_0^{T^n}b(\chi_s,\omega)ds\right]}{E^{\hat{P}}\big[T^n\big]} & = & 
E^{Q^s}\big[Z_1^s\big]
\end{aligned}
\end{equation}
holds and then we find that the sequence on the left is in fact constant and equals $v$.
Since the measure $\hat{P}\times K_0$ is absolutely continuous 
with respect to $P\times K_0$ it follows from (\ref{293}) by using (1) of Theorem \ref{47} and the fact that 
$P\times K_0$-a.s, $W_t/t\longrightarrow 0,$ as $t\to\infty$,
\begin{equation*}
\hat{P}\times K_0\mbox{-a.s.,}\qquad \frac{1}{t}\int_0^t b(\chi_s,\omega)ds\overset{\scriptscriptstyle{t\to\infty}}
{\longrightarrow} E^{Q^s}\big[Z_1^s\big].
\end{equation*}
By dominated convergence this limit holds true in $L^1(\hat{P}\times K_0)$ as well. 
Because of the ergodicity of $(W{{\cap}}\{0\in\cal{C}\},\hat{\theta}_1,\hat{P})$, which is a 
consequence of the ergodicity of $(W,\theta_1,P),$ see (34) on page 357 in \cite{N}, we have:  
\begin{equation}\label{68}
\frac{T^n}{n}\overset{\scriptscriptstyle{(\ref{0200})}}{=}\frac{1}{n}\sum_{k=0}^{n-1}T^1\circ\hat{\theta}_k
\overset{\scriptscriptstyle{n\to\infty}}{\longrightarrow} 
E^{\hat{P}}[T^1]\overset{\scriptscriptstyle{(\ref{21}),(\ref{24})}}
{<}\infty\qquad \hat{P}\mbox{-a.s. and in }L^1(\hat{P}),
\end{equation}
and we find that $\hat{P}\times K_0$-a.s. and in $L^1(\hat{P}\times K_0),$
\begin{equation*}
\lim_{n\to\infty}\frac{1}{n}\int_0^{T^n}b(\chi_s,\omega)ds=
\lim_{n\to\infty}\frac{T^n}{n}\ \frac{1}{T^n}\int_0^{T^n}b(\chi_s,\omega)ds=E^{Q^s}\big[Z_1^s\big]E^{\hat{P}}\big[T^1\big].
\end{equation*}Together with (\ref{0200}) and (\ref{23}), (\ref{65}) now follows. 
For a fixed $(w,\lambda_.)\in W\cap\{0\in\cal{C}\}$ and $k\geq 1,$ 
we find by an application of (\ref{257}) with $n=T^k$ and then with $n=T^k-1$ 
and similar considerations to those leading to (\ref{701}) that
\begin{eqnarray*}
&& E^{K_0}\left[\int_{T^k}^{T^{k+1}}b(\chi_u,\omega)du\right]=  E^{K_0}\left[\int_{0}^{T^{1}\circ\hat{\theta}_k}b(\chi_{T^k+u},\omega)du\right]\\
&&=\int_{\bb{R}^{d_2}}\frac{dy}{vol(d_2)}\bb{E}\bigg[
E^{K_{0,\omega}}\bigg[\bbm{1}_{\{y\in B^{d_2}_{1}(X^2_{T^k-1})\}}\bigg]
E^{K_{0,\bar{\omega}}\circ\hat{\theta}_k}
\bigg[\int_0^{T^1\circ\hat{\theta}_k}b\big((X^1_u\circ\hat{\theta}_k,X^2_u),\bar{\omega}\big)du\bigg]\bigg]
\end{eqnarray*}with $\bar{\omega}=\tau_{(w(T^k),y)}(\omega).$ By an independence argument as above (\ref{702}) and 
stationarity of the environment we finally obtain that\begin{equation}\label{260}
E^{K_0}\left[\int_{T^k}^{T^{k+1}}b(\chi_u,\omega)du\right]=
E^{K_0}\left[\int_{0}^{T^{1}}b(\chi_u,\omega)du\right]\circ\hat{\theta}_k.
\end{equation}
Recalling that the measure $\hat{P}$ is invariant under $\hat{\theta}_k,$ see (\ref{23}), we thus find\begin{eqnarray*}
E^{\hat{P}\times K_0}\left[\int_0^{T^n}b(\chi_u,\omega)du\right] & = & 
\sum_{k=0}^{n-1}E^{\hat{P}\times K_0}\left[\int_{T^k}^{T^{k+1}}b(\chi_u,\omega)du\right]\\
& \overset{(\ref{23}),(\ref{260})}{=} & n E^{\hat{P}\times K_0}\left[\int_0^{T^1}b(\chi_u,\omega)du\right],
\end{eqnarray*}and \begin{equation*}E^{\hat{P}}[T^n]\overset{(\ref{0200})}{=}E^{\hat{P}}\bigg[\sum_{k=0}^{n}
T^1\circ\hat{\theta}_k\bigg]\overset{(\ref{23})}{=}
nE^{\hat{P}}[T^1],\end{equation*} which shows that the sequence in (\ref{65}) is indeed constant and equal to 
\begin{equation}\label{340}
v\df\frac{E^{\hat{P}\times K_0}[\int_0^{T^1}b(\chi_u^{},\omega)du]}{E^{\hat{P}}[T^1]}=
\frac{E^{\hat{Q}^s}[\int_0^{T^1}b(Z^s_u,\omega)du]}{E^{\hat{P}}[T^1]}{=}
\frac{E^{\hat{Q}^s}[\int_0^{T^1}b(\alpha_u^s)du]}{E^{\hat{P}}[T^1]},
\end{equation}
where we used Remark \ref{540} in first equality and definition (\ref{45}) together with the fact that $Z^s_{T^0}=Z^s_0=0,$ 
$\hat{Q}^s$-a.s., 
in the second equality in (\ref{340}). The second and the third equality in (\ref{55}) then follows from 
(\ref{340}) by applying (\ref{53}) and (\ref{234}) to the last expression in (\ref{340}).

\begin{flushright}$\Box$\end{flushright}

{\rem{The formula for the limiting 
velocity, see (\ref{55}), is reasonably explicit and depends only on a finite piece of trajectory up to the first
cut time after time 0 and its first moment.}}

\section{Two invariance principles under the annealed measure}\label{c}

In  this section we provide two central limit theorems under the annealed measure. The first one is shown under a
symmetry assumption on the drift and $d_1\geq 7$, see Theorem \ref{103}, whereas for the second theorem there is no 
symmetry assumption but we need to assume that $d_1\geq 13.$\\

For integer $n\geq 1$ we denote with $I_{n}$ the $n\times n$-dimensional identity matrix. We further introduce the reflection
\begin{eqnarray*}
\cal{R}:\bb{R}^{d_1}\times\bb{R}^{d_2}&\longmapsto &\bb{R}^{d_1}\times\bb{R}^{d_2}\\
(x,y) & \longmapsto & (x,-y).
\end{eqnarray*}
For the first central limit theorem we assume the following antipodal symmetry in the last $d_2$ components 
of the drift under the measure $\bb{P}:$
\begin{equation}\label{100}
\left(\cal{R}(b(z,\omega))\right)_{z\in\bb{R}^{d}}\mbox{ has the same law as }\left(b(\cal{R}(z),\omega)\right)_{z\in\bb{R}^{d}}.
\end{equation}
Since the first $d_1$ components of the drift $b(\cdot,\cdot)$ vanish, we have that $\cal{R}\left(b(\cdot,\cdot)\right)$ equals $-b(\cdot,\cdot).$
 Note that when (\ref{100}) holds, then $\cal{R}(X_.)$ has the same law under $P_0$ as $X_.,$ and $E_0[X_t]=0$ for all $t\geq 0.$ 
By definition of $(W^{'}_t)_{t\geq 0},$ see below (\ref{56}), we have that $P_{0,\omega}\mbox{-a.s., }X_t=X_0+\int_0^t b(X_s,\omega)ds+W^{'}_t,
$ for each $\omega\in \Omega.$ 
The strong law of large numbers for Brownian motion, see Problem 9.3 in \cite{KS}, and Theorem \ref{89} imply that 
$P_0\mbox{-a.s., }\frac{1}{t}\int_0^t b\left(X_s,\omega\right)ds\longrightarrow v,$ 
and hence with dominated convergence the convergence holds in $L^1(P_0)$ as well. So $E_0[X_t]=E_0[\int_0^t b(X_s,\omega)ds]=0,$
and we deduce that 
the limiting velocity in (\ref{55}) vanishes under the assumption (\ref{100}).\\ 
{\rem{\label{0150}
A possible example of a drift $b^*(z,\omega),$ with $z\in\bb{R}^d,\omega\in\Omega,$ see (\ref{77}), such that (\ref{100}) is satisfied
can be constructed as follows. We consider a canonical Poisson point process on $\bb{R}^d$ with constant intensity
 as the random environment. Pick an $\bb{R}^{d_2}$-valued 
measurable function $\varphi(z),z\in\bb{R}^d,$ which is supported in a ball 
of radius $R/4$ and such that $\varphi(\cal{R}(z))=-\varphi(z)$ holds for all $z\in\bb{R}^d.$
Then make the convolution of the Poisson point process 
with the function $\varphi$ and truncate the new function. After smoothing out with a 
Lipschitz continuous real-valued mollifier $\rho(z), z\in\bb{R}^d,$
supported in a ball of radius $R/4$ and such that $\rho(\cal{R}(z))=\rho(z)$ for all $z\in\bb{R}^d,$ one obtains an 
example of a possible $b^*(z,\omega).$ 
}}\\

For two $C(\bb{R}_+,\bb{R}^d)$-valued sequences $\xi^n_.$ and 
$\zeta_.^n,\ n\geq 1,$ respectively defined on the probability spaces $(\Xi_1,\cal{D}_1,\mu_1)$ and  $(\Xi_2,\cal{D}_2,\mu_2)$  
we say that $(\xi^n_.)_{n\geq 1}$ under $\mu_1$ is weak convergence equivalent (abbreviated by wce) to  $(\zeta^n_.)_{n\geq 1}$ under 
$\mu_2,$ if the weak convergence of the law of $\xi_.^n$ under $\mu_1$ is equivalent to the weak convergence of the law of 
$\zeta_.^n$ under $\mu_2,$ and if both limits are the same, when weak convergence holds true.\\

Before we come to the main results of this section we briefly discuss some integrability properties stated in the 
following

{\lemma{\begin{eqnarray}
\label{032}&& \mbox{For all }\eta \geq 1:\ T^1\in L^{\eta}(P)\Leftrightarrow T^1\in L^{\eta+1}(\hat{P})\\
\label{111}&& T^1\in L^2(\hat{P}),\mbox{ when }d_1\geq 7.\\
\label{600}&& T^1\in L^4({P})\mbox{ and } T^1\in L^5(\hat{P}),\mbox{ when }d_1\geq 13.\\
\label{601}&& \sup_{{{a\in [0,T^1]}}}|\chi_{a}^{}|\in L^2(\hat{P}\times K_0),\mbox{ when }d_1\geq 7.\\
\label{602}&& \sup_{{{a\in [0,T^1]}}}|\chi_{a}^{}|\in L^4(\hat{P}\times K_0),\mbox{ when }d_1\geq 13.
\end{eqnarray}
}}\\

{\textbf{Proof:}} The equivalence (\ref{032}) is an easy consequence of (\ref{25}). With the help of (\ref{26}) we 
find that $T^1\in L^1{(P)}$ when $d_1\geq 7$ and  $T^1\in L^4{(P)}$ when $d_1\geq 13$ and so, (\ref{032}) yields 
(\ref{111}) and (\ref{600}). With the help of the integral representation of $\chi_.^{},$ see (\ref{150}) 
and note that $\hat{P}\times K_{0,\omega}\ll {P}\times K_{0,\omega},$ and (\ref{2}) we see that for each $\omega\in\Omega,\
 \hat{P}\times K_{0,\omega}$-a.s.,\begin{equation}\label{603}
 \sup_{{{a\in [0,T^1]}}}|\chi_{a}^{}|^2\leq 2\kappa^2 (T^1)^2+2\sup_{{{a\in [0,T^1]}}}|W_{a}|^2.
 \end{equation}Taking the $\hat{P}\times K_{0,\omega}$-expectation on both sides of (\ref{603}) we observe that (\ref{601}) 
 follows from (\ref{111}) and once we show that uniformly in $\omega,$\begin{equation}\label{604}
 E^{\hat{P}\times K_{0,\omega}}\left[\sup_{{{a\in [0,T^1]}}}|W_{a}|^2\right]\leq c_5(\varepsilon)<\infty.
 \end{equation}The left-hand side of (\ref{604}) is equal to\begin{equation}\label{033}
\sum_{n\geq 1}E^{\hat{P}\times K_{0,\omega}}
\left[\sup_{{\scriptscriptstyle{a\in [0,n]}}}|W_{a}|^2, T^1=n\right]
\overset{{\mbox{\tiny{H\"older}}}}{\leq}\sum_{n\geq 1}E^{\hat{P}\times K_{0,\omega}}
\left[\sup_{{\scriptscriptstyle{a\in [0,n]}}}|W_{a}|^{2p}\right]^{1/p}\hat{P}[T^1=n]^{1/q},
 \end{equation}with $1<q<\frac{6}{5}$ and $p$ the conjugate exponent. From (\ref{21}) and the definition of $\hat{P},$ see (\ref{23}), 
we see that\begin{equation}\label{620}
\hat{P}[\,\cdot\,]\leq c_1(\varepsilon)^{-1}P[\,\cdot\,].
\end{equation}An application of Burkholder-Davis-Gundy-Inequality, see p. 166 of \cite{KS}, yields
\begin{equation*}
 E^{\hat{P}\times K_{0,\omega}}\left[\sup_{{\scriptscriptstyle{a\in [0,n]}}}|W_{a}|^{2p}\right]^{1/p}\overset{\mbox{\tiny{(\ref{620})}}}{\leq}
 c(\varepsilon,q)E^{{P}\times K_{0,\omega}}\left[\sup_{{\scriptscriptstyle{a\in [0,n]}}}|W_{a}|^{2p}\right]^{1/p}\leq  c(\varepsilon,q)n\ ,
 \end{equation*}and hence the right-hand side of (\ref{033}) is less or equal to\begin{equation}\label{605}
c(\varepsilon,q)\sum_{n\geq 1}n\hat{P}[T^1=n]^{1/q}
= c(\varepsilon,q)\sum_{n\geq 1}n\hat{P}[T^1=n]^{1/2}\hat{P}[T^1=n]^{1/q-1/2}.
\end{equation}From an application of Cauchy-Schwarz' inequality follows that (\ref{605}) is dominated 
by\begin{equation}\label{522}
c(\varepsilon,q)\bigg\{E^{\hat{P}}[(T^1)^2]^{1/2}\bigg(\sum_{n\geq 1}\hat{P}[T^1=n]^{2/q-1}\bigg)^{1/2}\bigg\}.
\end{equation}Since $\hat{P}[T^1=n]\leq c_1(\varepsilon)^{-1} P[T^1>n-1]$ holds and $d_1\geq 7,$ one easily checks by using 
(\ref{26}) that the sum in (\ref{522}) with $1<q<\frac{6}{5}$ is bounded by a constant $c(\varepsilon,q)$ and with 
(\ref{111}), (\ref{604}) then follows. (\ref{602}) is shown analogously to (\ref{601}) with $1<q<\frac{6}{5}$ 
using now (\ref{600}) instead of (\ref{111}).\begin{flushright}$\Box$\end{flushright}

We now are ready to state our first invariance principle.

{\thm{\label{103}Let us assume $d_1\geq 7$ and (\ref{100}). Under the measure $P_0$, the $C(\bb{R}_+,\bb{R}^{d})$-valued 
random variables
\begin{equation}\label{900}
B^r_.\df \frac{1}{\sqrt{r}}X_{r\cdot},\quad r>0,
\end{equation}
converge in law to a $d$-dimensional Brownian motion $B_.$ with covariance matrix
\begin{equation}\label{101}
{A}=E^{\hat{P}}[T^1]^{-1}\left(\begin{array}{cc}
E^{\hat{P}}[T^1]I_{d_1} & 0 \\
0 & E^{\hat{Q}^s}[(Y^s_{T^1})(Y^s_{T^1})^t]
\end{array}\right)\in \bb{R}^{d\times d},
\end{equation} as $r\to\infty.$}}\\

{\rem{\label{302}Before giving the proof of the theorem, let us recall some classical facts about weak convergence on 
$C(\bb{R}_+,\bb{R}^d),$ that will be used several times throughout Section \ref{c}. More details on the following results 
can be found in Chapter 3 of \cite{EK} and in Section 3.1 of \cite{S}. Let us consider the space $C(\bb{R}_+,\bb{R}^{d})$ 
and the metric
\begin{equation*}
d(\xi_.,\zeta_.)\df\sum_{m=1}^{\infty}2^{-m}\sup_{0\leq t\leq m}(|\xi_t-\zeta_t|\wedge 1)\leq 1,\quad \xi_.,\zeta_.\in 
C(\bb{R}_+,\bb{R}^d).
\end{equation*}
Then $C(\bb{R}_+,\bb{R}^d)$ with the topology induced by $d(\cdot,\cdot)$ is a Polish space. Suppose $\xi_.^n$ and $\zeta_.^n, n\geq 1,$ are 
two $C(\bb{R}_+,\bb{R}^d)$-valued sequences on some probability space $(\Xi,\cal{D},\mu).$ If $d(\xi_.^n,\zeta_.^n)$ 
converges in $\mu$-probability to 0, then $(\xi^n_.)_{n\geq 1}$ under $\mu$ is wce to $(\zeta^n_.)_{n\geq 1}$ under 
$\mu,$ see below Remark \ref{0150} for the meaning of wce. Note that in order to verify the convergence in probability $\mu$  
of the distance $d(\xi_.^n,\zeta_.^n)$ to 0, it suffices to check that for any $T>0,\varepsilon>0,$
\begin{equation}\label{102}
\mu\left(\sup_{0\leq t\leq T}|\xi_t^n-\zeta_t^n|>\varepsilon\right)\overset{\scriptscriptstyle{n\to\infty}}
{\longrightarrow}0.
\end{equation}}}

{\textbf{Proof of Theorem \ref{103}}:} 
Observe that Theorem \ref{103} follows if we show that for $n\geq 1$ integer,
\begin{equation}\label{301}B^n_.\longrightarrow B_.\mbox{ in law under }P_0,\mbox{ as }n\to\infty.
\end{equation}
Indeed, (\ref{301}) implies that for $s_n\nearrow\infty,$ the sequence $\ 
[s_n]^{-1/2}X_{[s_n]\cdot}$ and thus $s_n^{-1/2}X_{[s_n]\cdot}$ converges in law to $B_.,$ recall (\ref{300}). 
Therefore, the laws of $s_n^{-1/2}X_{[s_n]\cdot}$ are tight and hence, by Theorem 2.4.10 of \cite{KS}, for all 
$T>0,\varepsilon>0,$ there exists an $\eta>0$
such that 
\begin{equation*}
\sup_{n\geq 1}P_0\left[\sup_{\substack{\scriptscriptstyle{|s-t|\leq\eta}\\ 
\scriptscriptstyle{0\leq s,t\leq T}}}\frac{1}{\sqrt{s_n}}|X_{[s_n]t}-X_{[s_n]s}|\geq\varepsilon\right]
\leq\varepsilon.
\end{equation*}
Since $\sup_{t\leq T}|t-\frac{s_n}{[s_n]}t|\overset{\scriptscriptstyle n\to\infty}{\longrightarrow}0,$ we obtain that for large $n,$
\begin{equation*}
P_0\left[\sup_{0\leq t\leq T}\frac{1}{\sqrt{s_n}}|X_{[s_n]t}-X_{s_nt}|\geq\varepsilon\right]\leq \varepsilon.
\end{equation*}
In view of Remark \ref{302}, this shows that $B_.^{s_n}$ converges in law to $B_.$ for any $s_n\nearrow\infty,$ which 
proves Theorem \ref{103}. For integer $n\geq 1$ we introduce the following piece-wise linear processes 
(recall the definitions (\ref{40}) and (\ref{43})):\begin{equation}\label{09}\begin{array}{rcl}
\bar{B}_.^n & \df & \frac{1}{\sqrt{n}}\Big\{X_{[n\cdot]}+(n\cdot-[n\cdot])\left(X_{[n\cdot]+1}^{}
-X_{[n\cdot]}^{}\right)\Big\},\\
\bar{Z}_n(\cdot) & \df & Z_{[n\cdot]}+(n\cdot-[n\cdot])\left(Z_{[n\cdot]+1}^{}-Z_{[n\cdot]}^{}\right),\\
\bar{Z}_n^s(\cdot) & \df & Z^s_{[n\cdot]}+(n\cdot-[n\cdot])\left(Z^s_{[n\cdot]+1}-Z^s_{[n\cdot]}\right).
\end{array}\end{equation}Note that the processes $\bar{B}^n_.,\frac{1}{\sqrt{n}}\bar{Z}_n(\cdot)$ and $\frac{1}{\sqrt{n}}\bar{Z}_n^s(\cdot)$ are the 
polygonal interpolations of $(X_k)_{k\geq 0},\ (Z_k)_{k\geq 0}$ and $(Z_k^s)_{k\geq 0}$ respectively, 
which are then rescaled in time and space as in the definition of $B^n_.$ for integers $n\geq 1,$ see (\ref{900}).
{\lemma{\label{510}$(B^n_.)_{n\geq 1}$ under $P_0$ is wce to $\left(\frac{1}{\sqrt{n}}\bar{Z}^s_n(\cdot)\right)
_{n\geq 1}$ under $Q^s.$}}\\

{\textsc{Proof:}} As a first step we show that\begin{equation}\label{321}
(B^n_.)_{n\geq 1}\mbox{ under }P_0\mbox{ is wce to }(\bar{B}^n_.)_{n\geq 1}\mbox{ under }P_0.
\end{equation} In view of Remark \ref{302}, see in particular (\ref{102}), it suffices to prove that for any $T>0$ the sequence of random variables
$\sup_{0\leq t\leq T}|B^n_t-\bar{B}^n_t|$ converges in $P_0$-probability to 0, as $n\to\infty.$ Indeed, the process $(W'_t)_{t\geq 0}$ 
defined below (\ref{56}) is a $d$-dimensional Brownian motion 
under $P_{0,\omega},$ and so for $T>0$ and $\varepsilon>0,$ when $n$ is large uniformly in $\omega,$
\begin{eqnarray*}
P_{0,\omega}\left[\sup_{0\leq t\leq T}|B^n_t- \bar{B}^n_t|\geq 4\varepsilon\right]&\leq &
P_{0,\omega}\left[\sup_{\substack{\scriptscriptstyle{k=0,\ldots,[Tn]}\\ 
\scriptscriptstyle{0\leq a\leq 1}}}|X_{k+a}-X_{k}|\geq 2\varepsilon\sqrt{n}\right]\\
&=& P_{0,\omega}\left[\sup_{\substack{\scriptscriptstyle{k=0,\ldots,[Tn]}\\ 
\scriptscriptstyle{0\leq a\leq 1}}}
  \bigg|\int_k^{k+a}b(X_s,\omega)ds+W^{'}_{k+a}-W^{'}_{k}\bigg|\geq 2\varepsilon\sqrt{n}\right]\\
&\leq & c(1+Tn)\exp\{-\frac{\varepsilon^2}{2d^2}n\},\end{eqnarray*}
where we used (\ref{2}) and Bernstein's inequality in the last line, and (\ref{321}) follows. From the identities
 in law stated in (1) of Theorem \ref{47} and Proposition \ref{59} we immediately deduce that
\begin{equation}\label{106}
(\bar{B}^n_.)_{n\geq 1}\mbox{ under }P_0\mbox{ is identical in law to }\left(\frac{1}{\sqrt{n}}\bar{Z}_n(\cdot)\right)_{n\geq 1}
\mbox{ under }Q^0.
\end{equation} A combination of (\ref{321}) and (\ref{106}) yields Lemma \ref{510} once we show that
\begin{equation}\label{810}
 \left(\frac{1}{\sqrt{n}}\bar{Z}_n(\cdot)\right)_{n\geq 1}\mbox{ under }Q^0\mbox{ is wce to }
 \left(\frac{1}{\sqrt{n}}\bar{Z}_n^s(\cdot)\right)_{n\geq 1}\mbox{ under }Q^s.
\end{equation}As in the proof of Theorem \ref{89} we define the processes $\bar{Z}_n(\cdot)$ and $\bar{Z}_n^s(\cdot)$ 
on a common probability space, see (\ref{60}) and below. Then we can again use the strategy discussed in Remark 
\ref{302} to prove (\ref{810}). In the notation (\ref{60})-(\ref{107}), using the fact that 
(\ref{63}) holds true, we find for $T>0$ that $Q$-a.s.,\begin{equation}\label{511}
\sup_{0\leq t\leq T}\frac{1}{\sqrt{n}}|\bar{Z}_n(t)\circ\pi^0-\bar{Z}_n^s(t)\circ\pi^s|\leq
\sup_{0\leq t\leq \frac{T^1}{n}}\frac{1}{\sqrt{n}}|\bar{Z}_n(t)\circ\pi^0-\bar{Z}_n^s(t)\circ\pi^s|.
\end{equation}Since $\bar{Z}_n(t)\circ\pi^0-\bar{Z}_n^s(t)\circ\pi^s,\ t\in [0,\frac{T^1}{n}],$ is a continuous 
process which is piece-wise linear between the times $0,\frac{1}{n},\frac{2}{n},\ldots,\frac{T^1}{n},$ we find 
that the right-hand side of (\ref{511}) is equal to \begin{equation*}
\sup_{k=0,\frac{1}{n},\ldots,\frac{T^1}{n}}\frac{1}{\sqrt{n}}|\bar{Z}_n(k)\circ\pi^0-\bar{Z}_n^s(k)\circ\pi^s|=
\sup_{k=0,\ldots,T^1}\frac{1}{\sqrt{n}}|{Z}(k)\circ\pi^0-{Z}^s(k)\circ\pi^s|,
\end{equation*}{which converges $Q$-a.s. to zero, as $n\to\infty$, since $Q[T^1<\infty]=1$, see (\ref{290}) and 
(\ref{22}). This concludes the proof of (\ref{810}) and thus of Lemma \ref{510}.\begin{flushright}$\Box$\end{flushright}}

\pagebreak

Let us define an integer-valued function $0\leq \varphi(t)$ tending to infinity $P$-a.s., such that
\begin{equation}\label{108}
T^{\varphi(t)}\leq t<T^{\varphi(t)+1}\quad\mbox{for all }t\geq 0,
\end{equation}and
\begin{equation}\label{109}
\Sigma_m\df Z^s_{T^m}-Z^s_{T^0},\quad m\geq 0.
\end{equation} Furthermore, let us introduce the polygonal interpolation of $\Sigma_m,m\geq 0:$\begin{equation}\label{02}
\bar{\Sigma}_{\cdot}\df\Sigma_{[\cdot]}+(\cdot-[\cdot])\left(\Sigma_{[\cdot]+1}-\Sigma_{[\cdot]}\right),
\end{equation}and for integer $n\geq 1,$
\begin{equation}\label{800}
\bar{\Sigma}^{\varphi}_{n}(\cdot)\df\Sigma_{\varphi(n\cdot)}+(n\cdot-[n\cdot])\left(\Sigma_{\varphi(n\cdot+1)}-\Sigma_{\varphi(n\cdot)}\right),
\end{equation}which is constant and equal to $\Sigma_{\varphi(T^k)}=\Sigma_k$ on the time interval 
$[\frac{T^k}{n},\frac{T^{k+1}}{n}-\frac{1}{n}),\ k\geq 0,$ 
and linear on the interval $[\frac{T^{k+1}}{n}-\frac{1}{n},\frac{T^{k+1}}{n}),$ interpolating the points $\Sigma_k$ and $\Sigma_{k+1}.$
{\lemma{\label{530}$\left(\frac{1}{\sqrt{n}}\bar{Z}_n^s(\cdot)\right)_{n\geq 1}$ under $Q^s$ is wce to
$\left(\frac{1}{\sqrt{n}}\bar{\Sigma}_n^{\varphi}(\cdot)\right)_{n\geq 1}$ under $Q^s.$}}\\

{\textsc{Proof:}} In view of Remark \ref{302} it suffices to show that for 
any $T>0$ and $\varepsilon>0$ the following probability converges to 0, as $n\to\infty:$
\begin{equation}\label{535}
Q^s\left[\sup_{0\leq t\leq T}\frac{1}{\sqrt{n}}|\bar{Z}^s_n(t)-\bar{\Sigma}^{\varphi}_n(t)|>4\varepsilon\right]
\leq Q^s\Bigg[\underbrace{\sup_{\substack{\scriptscriptstyle{k=0,\ldots,[Tn]+1}\\
\scriptscriptstyle{a\in [0,T^{k+1}-T^k]}}}|Z^s_{T^k+a}-Z^s_{T^k}|>
\varepsilon{\sqrt{n}}}_
{\df A_n}\Bigg].
\end{equation}
Since the event $A_n$ is invariant under the shift $\Theta_{T^0}$ and the image of $Q^s$ under $\Theta_{T^0}$ is $E^{\hat{Q}^s}
[\, \cdot\, ,T^1]/E^{\hat{P}}[T^1],$ see (\ref{234}), it follows with the help of Cauchy-Schwarz' inequality that
\begin{equation}\label{128}
Q^s[A_n]=Q^s[\Theta_{T^0}^{-1}(A_n)]\leq E^{\hat{P}}[(T^1)^2]^{1/2}\hat{Q}^s[A_n]^{1/2}/E^{\hat{P}}[T^1],
\end{equation}where $E^{\hat{P}}[(T^1)^2]<\infty,$ see (\ref{111}). Thus, Lemma \ref{530} will follow once we show that
\begin{equation}\label{113}
\lim_{n\to\infty}\hat{Q}^s[A_n]=0.
\end{equation}Using (\ref{141}) and the fact that $\hat{\Theta}_k$ preserves $\hat{Q}^s,$ see the proof of Proposition \ref{88},
we find that\begin{equation}\label{533}
\hat{Q}^s[A_n]\leq (2+Tn)\hat{Q}^s\left[\sup_{a\in [0,T^1]}|Z_a^s|>\varepsilon\sqrt{n}\right]\leq\frac{2+Tn}{\varepsilon^2 n}
E^{\hat{Q}^s}\left[\sup_{a\in [0,T^1]}|Z_a^s|^2,\sup_{a\in [0,T^1]}|Z_a^s|>\varepsilon\sqrt{n}\right].
\end{equation}From (\ref{601}) and Remark \ref{540} follows that the last expression vanishes, as $n\to\infty,$ and hence (\ref{113}) 
holds true. This finishes the proof of Lemma \ref{530}.\begin{flushright}$\Box$\end{flushright}

\pagebreak

{\lemma{\label{01}
Under $Q^s,\ \left(\frac{1}{\sqrt{n}}\bar{\Sigma}_{n\cdot}\right)_{n\geq 1}$ converges in law to 
$\sqrt{E^{\hat{P}}\left[T^1\right]}B_.,$ as $n\to\infty.$ }}
Before proving Lemma \ref{01}, let us explain how we conclude the proof of Theorem \ref{103}. Once we show that
\begin{equation}\label{03}\left(\frac{1}{\sqrt{n}}\bar{\Sigma}_n^{\varphi}(\cdot)\right)_{n\geq 1}
\mbox{ under $Q^s$ is wce to }
\left(\frac{1}{\sqrt{n}}\bar{\Sigma}_{n\cdot/E^{\hat{P}}[T^1]}\right)_{n\geq 1}\mbox{ under }Q^s,
\end{equation}we find with Lemma \ref{01} and a transformation of time that the first sequence in (\ref{03}) converges 
weakly to $B_.,$ as $n\to\infty,$ and hence with Lemma \ref{530} and \ref{510} we deduce that (\ref{301}) holds, which finishes 
the proof of Theorem \ref{103}.\\
For the proof of (\ref{03}) first note that since $T^n-T^0$ is invariant under 
$\Theta_{T^0}$ we find by similar arguments to those leading to (\ref{128}) 
that the convergence in (\ref{68}) holds true $Q^s$-a.s. and not only $\hat{P}$-a.s.. It follows that $\frac{\varphi(t)}{t}\longrightarrow 
E^{\hat{P}}[T^1]^{-1},\ Q^s$-a.s. and hence with Lemma 9.2 on page 572 of \cite{G}:\begin{equation}\label{114}
\mbox{for all }T\geq 0,\ Q^s\mbox{-a.s., }\sup_{0\leq t\leq T}\bigg|\frac{\varphi(nt)}{n}-\frac{t}{E^{\hat{P}}[T^1]}
\bigg|\overset{\scriptscriptstyle{n\to\infty}}{\longrightarrow}0,
\end{equation}and so for $\varepsilon>0,\eta>0,T>0$ and $n$ large enough,\begin{equation}\label{04}
Q^s\left[\sup_{0\leq t\leq T}\left|\frac{\varphi(nt)}{n}-\frac{t}{E^{\hat{P}}[T^1]}\right|\geq \eta\right]\leq \varepsilon.
\end{equation}Furthermore, from Lemma \ref{01} we infer that the laws of $n^{-1/2}\bar{\Sigma}_{n\cdot}$ under 
$Q^s$ are tight and hence for all $T>0,\varepsilon >0,$ there exists an $\eta>0$ such that
\begin{equation}\label{05}
\sup_{n\geq 1}Q^s\left[\sup_{\substack{|s-t|\leq \eta\\ 0\leq s,t\leq T}}\frac{1}{\sqrt{n}}
\left|\bar{\Sigma}_{nt}-\bar{\Sigma}_{ns}\right|\geq \varepsilon\right]\leq\varepsilon,
\end{equation}see Theorem 2.4.10 of \cite{KS}. Together with (\ref{04}) we thus obtain that for arbitrary $\varepsilon>0$ and $ T>0,
$\begin{equation}\label{06}
Q^s\left[\sup_{ 0\leq t\leq T}\frac{1}{\sqrt{n}}
\left|\bar{\Sigma}_{nt/E^{\hat{P}}[T^1]}-\bar{\Sigma}_{\varphi(nt)}\right|\geq \varepsilon\right]\leq2\varepsilon
\end{equation}for sufficiently large $n.$ In order to prove (\ref{03}) it suffices to show that for $T>0$ and $\varepsilon>0$ the following probability
tends to zero with $n,$ see Remark \ref{302}: 
\begin{eqnarray*}
Q^s\left[\sup_{0\leq t\leq T}\Big|\bar{\Sigma}_{{nt}/{E^{\hat{P}}[T^1]}}-\bar{\Sigma}^{\varphi}_n(t)\Big|>2\varepsilon\sqrt{n}\right]
&\leq& Q^s\left[\sup_{0\leq t\leq T}\Big|\bar{\Sigma}_{{nt}/{E^{\hat{P}}[T^1]}}-\bar{\Sigma}_{\varphi(nt)}\Big|>\varepsilon\sqrt{n}\right]\\
&& +Q^s\left[\sup_{0\leq t\leq T}\Big|\bar{\Sigma}_{\varphi(nt)}-\bar{\Sigma}^{\varphi}_n(t)\Big|>\varepsilon\sqrt{n}\right].
\end{eqnarray*}The first expression on the right-hand side vanishes, as $n\to\infty,$ due to (\ref{06}). Moreover, it can easily be seen 
that the second expression on the right-hand side of the above inequality is less or equal to ${Q}^s[A_n],$ see 
(\ref{535}) for the definition of $A_n,$ which tends to 0, as 
$n\to\infty$, see (\ref{128}) and (\ref{113}). This finishes the proof of (\ref{03}) and hence of Theorem \ref{103}.\begin{flushright}$\Box$\end{flushright}

{\textsc{Proof of Lemma \ref{01}:}} Let us denote with $\bar{\Sigma}^{d_1}_.$ and $\bar{\Sigma}^{d_2}_.$ the first $d_1$ respectively 
the last $d_2$ components of the process $\bar{\Sigma}_.$ defined in (\ref{02}). Note that Lemma \ref{01} follows from the next two statements:
\begin{equation}\label{118}\begin{array}{l}
\mbox{under }P,\mbox{ the sequence } \frac{1}{\sqrt{n}}\bar{\Sigma}^{d_1}_{n\cdot}, n\geq 1,
\mbox{ converges in law to a $d_1$-dimensional}\\
\mbox{Brownian motion with covariance matrix }E^{\hat{P}}[T^1]I_{d_1},\mbox{ as }n\to\infty,
\end{array}\end{equation}and  for ${P}$-a.e. $(w,\lambda_.)\in W,$
\begin{equation}\label{116}\begin{array}{l}
\mbox{under the measure }M^s,\mbox{ the sequence }\ \frac{1}{\sqrt{n}}\bar{\Sigma}^{d_2}_{n\cdot},n\geq 1,
\mbox{ converges in}\\ \mbox{law to a }d_2\, \mbox{-dimensional Brownian motion with covariance matrix}\\
E^{\hat{Q}^s}[(Y^s_{T^1})(Y^s_{T^1})^t]\in\bb{R}^{d_2\times d_2}\mbox{ (independent of }(w,\lambda_.)),\mbox{ as }n\to\infty.\hspace{2.2cm}
\end{array}\end{equation}
Indeed, from (\ref{118}) and (\ref{116}) we can easily deduce that under $Q^s,$ the laws of 
$n^{-1/2}\bar{\Sigma}_{n\cdot},$ $n\geq 1,$  are tight,
 see Theorem 2.4.7 and 2.4.10 of \cite{KS}. Therefore, in order to prove Lemma \ref{01}, it suffices to show weak convergence of all finite dimensional 
distributions of $n^{-1/2}\bar{\Sigma}_{n\cdot}$ to the finite dimensional distributions of 
${\scriptstyle{\sqrt{E^{\hat{P}}[T^1]}}}B_.,$ as $n\to\infty,$ see Theorem 2.4.15 in \cite{KS}. But this can easily be inferred
from (\ref{118}) and (\ref{116}) with the help of characteristic functions.\\
Now, let us explain how to see (\ref{118}). Similarly as in the 
proof of (\ref{03}), we first note that for $\varepsilon>0,\eta>0,T>0$ and $n$ large enough,\begin{equation}\label{07}
P\left[\sup_{0\leq t\leq T}\left|\frac{T^{[nt]}}{E^{\hat{P}}[T^1]n}-t\right|\geq \eta\right]\leq \varepsilon.
\end{equation}Furthermore, by definition and self-similarity of Brownian motion we know that under $P,$ the processes 
$\ n^{-1/2}X^1_{{\scriptscriptstyle{E^{\hat{P}}[T^1]}}n\cdot}, n\geq 1,$ are distributed as a $d_1$-dimensional Brownian motion with covariance matrix 
$E^{\hat{P}}[T^1]I_{d_1},$ and hence their laws are tight. So, we can derive the same estimate as in (\ref{05}) but for the process 
$X^1_{{\scriptscriptstyle{E^{\hat{P}}[T^1]}}n\cdot}.$ Together with 
(\ref{07}) we thus obtain that\begin{equation}\label{08}
P\left[\sup_{0\leq t\leq T}\left|X^1_{T^{[nt]}}-X^1_{E^{\hat{P}}[T^1]nt}\right|\geq \varepsilon\sqrt{n}\right]\leq 2\varepsilon
\end{equation}for sufficiently large $n.$ Pick $T>0,\ \varepsilon>0$ and then observe that\begin{eqnarray*}
&&P\left[\sup_{0\leq t\leq T}\left|\bar{\Sigma}_{nt}^{d_1}-
X^1_{E^{\hat{P}}[T^1]nt}\right|> 3\varepsilon\sqrt{n}\right]\\
&&\leq P\left[\sup_{0\leq t\leq T}\left|X^1_{T^{[nt]}}-X^1_{E^{\hat{P}}[T^1]nt}\right|> \varepsilon\sqrt{n}\right] + 
 P\Bigg[\sup_{\substack{\scriptscriptstyle{k=0,\ldots,[Tn]}\\
\scriptscriptstyle{a\in [0,T^{k+1}-T^k]}}}\left|X^1_{T^{k}+a}-X^1_{T^{k}}\right|> \varepsilon\sqrt{n}\Bigg]. 
\end{eqnarray*}The first term after the above inequality tends to zero with $n$ due to (\ref{08}) whereas the second term is dominated by 
$Q^s(A_n)$ and thus converges to zero as well, as $n\to\infty.$ In view of Remark \ref{302} this proves (\ref{118}).\\
We now prove (\ref{116}). Note that for $\hat{P}$-a.e. $(w,\lambda_.)\in W
{{\cap}}\{0\in \cal{C}\},$ under $M^s,$ the increments $Y^s_{T^n}-Y^s_{T^{n-1}},\ n\geq 1,$ are independent, 
see (\ref{39}) and (\ref{44}), with mean zero, which is a consequence of the symmetry assumption (\ref{100}). Indeed, for $\hat{P}$-a.e. $(w,\lambda_.)\in W
{\cap}\{0\in \cal{C}\},$\begin{equation}\label{832}
E^{M^s}[Y^s_{T^n}-Y^s_{T^{n-1}}]=\left(E^{K_0}[X_{T^1}^2]\right)\circ\hat{\theta}_{n-1}\overset{\scriptstyle{(\ref{190})}}{=}
\bb{E}\left[E^{K_{0,\omega}}[X_{T^1}^2] \right]\circ\hat{\theta}_{n-1}.
\end{equation} In Remark \ref{034} we have seen that we can write\begin{equation}\label{035}
E^{K_{0,\omega}}[X_{T^1}^2]=
\int_{\bb{R}^{d_2}}\cdots\int_{\bb{R}^{d_2}}dy_1\cdots dy_{T^1}\prod_{k=0}^{T^1-1}h\left(w(k+\cdot)-w(k),
\lambda_k,y_k,y_{k+1},\hat{\omega}_k\right)y_{T^1},
\end{equation}with $\hat{\omega}_k=\tau_{(w(k),0)}(\omega),\ k=0,\ldots, T^1-1$ and $y_0:=0.$
Let us denote with $h^{(\cal{R})}$ the analogue to $h,$ see (\ref{28}), defined via the transition density 
$p^{(\cal{R})}_{w,\omega}(1,\cdot,\cdot)$ 
for $\omega\in\Omega,\ w\in W^{d_1}_0,$ attached by (\ref{13}) and (\ref{255}) to the drift 
$-b^{*}((w(\cdot),-\,\cdot),\omega),$ 
see also (\ref{77}).
Since $(X^2_t)_{t\in [0,1]}$ under $\tilde{P}_{y_k,y_{k+1}}$ has the same law 
as  $(-X^2_t)_{t\in [0,1]}$ under $\tilde{P}_{-y_k,-y_{k+1}},$ see below (\ref{79}) for the definition of
$\tilde{P}_{\cdot,\cdot},$ 
we see that for $k=0,\ldots,T^1-1,$\begin{equation*}
p^{(\cal{R})}_{w,\omega}(1,-y_k,-y_{k+1})=p_{w,\omega}(1,y_k,y_{k+1})
\end{equation*}and hence
\begin{equation}\label{0101}
h^{(\cal{R})}\left(w(k+\cdot)-w(k),\lambda_k,-y_k,-y_{k+1},\hat{\omega}_k\right)=
h\left(w(k+\cdot)-w(k),\lambda_k,y_k,y_{k+1},\hat{\omega}_k\right).
\end{equation} With the help of Theorem 44 on page 158 in \cite{P} one can see that for fixed $(w,\lambda_.)\in W$ 
and $y_k,y_{k+1}\in\bb{R}^{d_2}$ the expression on the left-hand side of (\ref{0101}) is a measurable function of 
$\cal{R}(b(\cal{R}(\,\cdot\,) ,\omega)),$ see (\ref{77}). 
So, (\ref{0101}) together with our symmetry assumption (\ref{100}) 
implies that under $\bb{P},$ the product in (\ref{035}) is identical in law to
\begin{equation*}
\prod_{k=0}^{T^1-1}h\left(w(k+\cdot)-w(k),\lambda_k,-y_k,-y_{k+1},
\hat{\omega}_k\right)
\end{equation*}and so, a transformation of variables $(y_1,\ldots,y_{T^1})\mapsto (-y_1,\ldots,-y_{T^1})$ and 
Fubini's Theorem then yield 
$E^{K_{0}}[X_{T^1}^2]=-E^{K_{0}}[X_{T^1}^2]=0.$ Hence in view of (\ref{832}),  $\hat{P}$-a.e. $(w,\lambda_.)\in W
{{\cap}}\{0\in \cal{C}\},$
\begin{equation*}
E^{M^s}[Y^s_{T^n}-Y^s_{T^{n-1}}]=0,
\end{equation*}since $\hat{\theta}_{n-1}$ preserves $\hat{P}.$ Furthermore,\begin{equation}\label{831}\begin{array}{rcl}
E^{\hat{P}}\left[E^{M^s}\left[|Y^s_{T^n}-Y^s_{T^{n-1}}|^2\right]\right]&\overset{\scriptstyle{(\ref{37}),(\ref{43})}}{\leq}&
E^{\hat{Q}^s}\left[|Z^s_{T^n}-Z^s_{T^{n-1}}|^2\right]\\[11pt]
&\overset{\scriptstyle{(\ref{141})}}{=}&E^{\hat{Q}^s}\left[|Z^s_{T^1}|^2\circ\hat{\Theta}_{n-1}\right]<\infty,
\end{array}\end{equation}since $\hat{\Theta}_{n-1}$ preserves $\hat{Q}^s$ and because of the integrability property (\ref{601}) and Remark \ref{540}. 
In particular it follows that for $\hat{P}$-a.e. $(w,\lambda_.)\in W{{\cap}}\{0\in \cal{C}\},$\begin{equation}
\label{0100}
Y^s_{T^n}-Y^s_{T^{n-1}}\in L^2(M^s(w,\lambda_.)).
\end{equation}
Note that $(W{{\cap}}\{0\in\cal{C}\},\hat{\theta}_1,\hat{P})$ is ergodic as a consequence of the
ergodicity of $(W,{\theta}_1,{P}),$ see (34) on page 357 in \cite{N}.
An application of an invariance principle for vector-valued, square-integrable martingale 
differences, see Theorem \ref{080}, shows that for $\hat{P}$-a.e. 
$(w,\lambda_.)\in W{{\cap}}
\{0\in\cal{C}\},$ under the measure $M^s(w,\lambda_.),$ the $C(\bb{R}_+,\bb{R}^{d_2})$-valued random variables $n^{-1/2}{\bar{\Sigma}}_{n\cdot}^{d_2},\ n\geq 1$ 
converge weakly to a $d_2$-dimensional Brownian motion with covariance matrix as in (\ref{116}), as $n\to\infty.$  
Note that in fact under the measure $M^s(w,\lambda_.)$ the increments $Y^s_{T^n}-Y^s_{T^{n-1}},\ 
n\geq 1,$ are independent and hence the standard functional central limit theorem for independent increments, 
which is an immediate consequence of Theorem \ref{080}, could be applied. The ergodicity of $(W{{\cap}}\{0\in\cal{C}\},\hat{\theta}_1,\hat{P})$ and the integrability 
property (\ref{831}) are used to show that the conditions (\ref{081}) and (\ref{082}) are satisfied.
Since for a continuous, bounded function $f$ on $W^{d_2}_+,$ the random variable 
{$E^{M^s}[f({n^{-1/2}}\bar{\Sigma}_{n\cdot}^{d_2}))]$} is invariant under $\theta_{T^0}$ and since the image of 
$P$ under $\theta_{T^0}$ is absolutely continuous with respect to $\hat{P},$ see (\ref{25}), 
it follows that (\ref{116}) holds in fact for $P$-a.e. $(w,\lambda_.)\in W.$ This finishes the proof of Lemma \ref{01}.
\begin{flushright}$\Box$\end{flushright}

The next theorem shows us that our model also contains examples of diffusions in random environment with possibly ballistic
behavior when $d_1\geq 13,$ satisfying an invariance principle, recall (\ref{77})-(\ref{78}).

{\thm{\label{120} Let $d_1\geq 13$ and recall the definition of $v$ in (\ref{55}). Under the measure $P_0,$ 
the $C(\bb{R}_+,\bb{R}^{d})$-valued random variables
\begin{equation*}
B^r_.\df\frac{1}{\sqrt{r}}\left(X_{r\cdot}-vr\cdot\right),\quad r>0,
\end{equation*}converge in law to a $d$-dimensional Brownian motion $B_.$ with covariance matrix $A$ given 
in (\ref{125}), as $r\to\infty.$}}\\

{\textbf{Proof:}} As in the proof of Theorem \ref{103} we can show that it suffices to prove that for 
$n\geq 1$ integer, \begin{equation}\label{801}
B^n_.\longrightarrow B_. \mbox{ in law under }P_0,\mbox{ as }n\to\infty,
\end{equation}see (\ref{301}). By similar arguments as in the proof of Lemma \ref{510} 
we find that \begin{equation}\label{060}\begin{array}{l}
\left(\bar{B}_.^n-v\sqrt{n}\,\cdot\right)_{n\geq 1}\mbox{ under }P_0
\mbox{ is wce to }\left(\frac{1}{\sqrt{n}}\left(\bar{Z}^s_n(\cdot)-vn\cdot\right)\right)_{n\geq 1}\mbox{ under }Q^s,
\end{array}\end{equation}see (\ref{106}) and (\ref{810}). Together with (\ref{321}) we see that (\ref{801}) 
follows once we show that\begin{equation}\label{122}\begin{array}{l}
\frac{1}{\sqrt{n}}\left(\bar{Z}^s_n(\cdot)-vn\cdot\right)\longrightarrow B_.\mbox{ in law under }Q^s,\mbox{ as }n\to\infty,
\end{array}\end{equation}
see (\ref{09}) for the definition of $\bar{Z}^s_n(\cdot)$ and $\bar{B}_.^n.$ In the notation\begin{equation}\label{121}
\cal{Z}\df Z^s_1-E^{Q^s}\left[Z^s_1\right]=Z^s_1-v,\quad\mbox{ we have that }\quad{Z}^s_n-vn
\overset{\scriptstyle{(\ref{52})}}{=}\sum_{k=0}^{n-1}\cal{Z}\circ\Theta_k,
\end{equation}recall (\ref{55}). We know from (\ref{252}) that\begin{equation}\label{040}
\cal{Z}\in L^m(Q^s)\qquad\mbox{ for all } m\in [1,\infty).
\end{equation}
For integers $k\geq 0,$ on the space $\Gamma^s$ we introduce the filtration\begin{equation}\label{123}
\cal{G}_k\df\sigma\left(Z^s_{n+1}-Z^s_n,\mbox{ for all } n\in\bb{Z}\mbox{ with }n<k\right).
\end{equation}The identity (\ref{141}) implies that for $k\geq 0:$\begin{equation}\label{084}
f\mbox{ is }\cal{G}_0\mbox{-measurable}\Longleftrightarrow f\circ\Theta_k\mbox{ is }\cal{G}_k\mbox{-measurable,}
\end{equation}and thus by stationarity, see (\ref{54}), we have that for $g\in L^1(Q^s),$\begin{equation}\label{092}
Q^s\mbox{-a.s. }E^{Q^s}\left[g\circ\Theta_k\,|\,\cal{G}_k\right]=E^{Q^s}\left[g\,|\,\cal{G}_0\right]\circ\Theta_k.
\end{equation}
The following adaptation of Gordin's method will play the key role in the proof.
{\lemma{\label{124}There is a $G\in L^2(\Gamma^s,\cal{G}_0,Q^s)$ such that\begin{equation*}
M_n\df G\circ\Theta_n-G+Z_n^s-vn=\sum_{k=0}^{n-1}\left(G\circ\Theta_1-G+\cal{Z}\right)\circ\Theta_k,\ n\geq 0,
\mbox{ is a }(\cal{G}_n)\mbox{-martingale.}
\end{equation*}}}Before we prove Lemma \ref{124}, let us explain how we conclude the proof of Theorem \ref{120} from it. 
Using stationarity of $\Theta_1$ under $Q^s,$ see (\ref{54}), and applying Chebychev's inequality we obtain 
for $T>0,\varepsilon>0,$
\begin{equation*}
Q^s\left[\sup_{\scriptscriptstyle{k=0,\ldots,[Tn]+1}}|G\circ\Theta_k|>\varepsilon\sqrt{n}\right]\leq
\frac{Tn+2}{\varepsilon^2 n}E^{Q^s}\left[|G|^2,|G|>\varepsilon\sqrt{n}\right]
\overset{\scriptscriptstyle{n\to\infty}}{\longrightarrow}0,
\end{equation*}
so that, in view of Remark \ref{302}, one easily finds that 
$n^{-1/2}\left(\bar{Z}^s_n(\cdot)-vn\cdot\right)_{n\geq 1}$ under $Q^s$ is wce to the rescaled polygonal 
interpolation of the process $M_k, k\geq 1,$ defined analogously to $\bar{B}_.^n$ in (\ref{09}), under $Q^s$ .
Since $M_n$ is a martingale with ergodic, square integrable increments, it follows from Theorem \ref{140}, see Appendix, 
that under the measure $Q^s,$ the rescaled polygonal interpolation of $M_k, k\geq 1,$ converges in law to a $d$-dimensional Brownian motion 
with covariance matrix\begin{equation}\label{125}
A=E^{Q^s}\left[\left(G\circ\Theta_1-G+\cal{Z}\right)\left(G\circ\Theta_1-G+\cal{Z}\right)^t\right],
\end{equation}
as $n\to\infty.$ This concludes the proof of (\ref{122}).
\begin{flushright}$\Box$\end{flushright}

{\textsc{Proof of Lemma \ref{124}:} First we explain how our 
claim follows once we show that\begin{equation}\label{126}
\sum_{k\geq 0}\big{\|}E^{Q^s}\left[\left(H\bbm{1}_{\{0\in\cal{C}\}}\right)\circ\Theta_k\ \big |\ \cal{G}_0\right]\big{\|}_{2}
<\infty,
\end{equation}
where in the previous notation, see (\ref{121}),\begin{equation}\label{127}
H\df\sum_{k=0}^{T^1-1}\cal{Z}\circ\Theta_k=\sum_{k=0}^{T^1-1}Z_1^s\circ\Theta_k-vT^1\overset{(\ref{52})}{=}
Z_{T^1}^s-vT^1.
\end{equation}
Note that $H\in L^2(Q^s).$ Indeed, \begin{eqnarray*}
E^{Q^s}\left[|H|^2\right]&\leq &E^{Q^s}\left[(T^1)^2\sum_{k=0}^{T^1-1}|\cal{Z}\circ\Theta_k|^2\right]=
\sum_{n\geq 1}n^2\sum_{k=0}^{n-1}E^{Q^s}\left[|\cal{Z}\circ\Theta_k|^2, T^1=n\right]\\
&\overset{{{\mbox{{\tiny{H\"older}}}}}}{\leq} & 
\sum_{n\geq 1}n^2\sum_{k=0}^{n-1}E^{Q^s}\left[|\cal{Z}\circ\Theta_k|^{2p}\right]^{1/p}
P[T^1=n]^{1/q},\end{eqnarray*}with $1<q<9/8$ and $p$ the conjugate exponent. Since $P[T^1=n]\leq 
P[T^1> n-1]$ and $E^{Q^s}[|\cal{Z}\circ\Theta_k|^{2p}]\leq c(p)<\infty$ by (\ref{54}) and (\ref{040}), we conclude 
with the help of (\ref{26})
that the right-hand side of the above inequality is finite when $d_{1}\geq 13.$\\ 
For $m\geq 1$ we define
\begin{equation}\label{085}
G^m\df E^{Q^s}\left[H\ |\ \cal{G}_0\right]+\sum_{k=1}^{m-1}E^{Q^s}\left[\left(H\bbm{1}_{\{0\in\cal{C}\}}\right)
\circ\Theta_k\ |\ \cal{G}_0\right].
\end{equation}
Then $G^m$ converges in $L^2(Q^s)$ towards a $G\in L^2(\Gamma^s,\cal{G}_0,Q^s)$ because of (\ref{126}). Moreover, for 
$m\geq 1$ we define $N_m=N((w,\lambda_.);[1,m-1])+1$ in the notation of (\ref{20}),  so that\begin{equation}\label{330}
\sum_{k=0}^{T^{N_m}-1}\cal{Z}\circ\Theta_k=H+\sum_{k=1}^{m-1}(H\bbm{1}_{\{0\in\cal{C}\}})\circ\Theta_k.
\end{equation}By stationarity, see (\ref{54}), we find that for $n\geq 0,$
\begin{equation}\label{129}
G\circ\Theta_n=\lim_{m\to\infty}G^m\circ\Theta_n\overset{\scriptstyle{(\ref{092}),(\ref{330})}}{=}
\lim_{m\to\infty}\ E^{Q^s}\left[\left(\sum_{k=0}^{T^{N_m}-1}\cal{Z}\circ\Theta_k\right)\circ\Theta_n
\ \Bigg|\ \cal{G}_n\right],
\end{equation}where the above limits are in $L^2(\Gamma^s,\cal{G}_n,Q^s).$
This yields for $n\geq 1,$\\[11pt]
$E^{Q^s}\left[M_{n+1}-M_n\ |\ \cal{G}_n\right]$
\begin{eqnarray*} 
&&=\lim_{m\to\infty}\ E^{Q^s}\left[\left(\sum_{k=0}^{T^{N_m}-1}\cal{Z}\circ\Theta_k\right)\circ\Theta_{n+1}+\cal{Z}\circ\Theta_n
-\left(\sum_{k=0}^{T^{N_m}-1}\cal{Z}\circ\Theta_k\right)\circ\Theta_{n}\ \Bigg|\ \cal{G}_n\right],
\end{eqnarray*}where the limit is in  $L^2(\Gamma^s,\cal{G}_n,Q^s).$ With the observation that\begin{equation*}
T^{N_m}\circ\theta_1=\left\{\begin{array}{ll}
T^{N_m}-1,& \mbox{ on }\{m\notin\cal{C}\},\\
T^{N_m+1}-1,& \mbox{ on }\{m\in\cal{C}\},
\end{array}\right.
\end{equation*}
we find that the quantity under the conditional expectation is equal to
\begin{equation*}
\left(\sum_{k=0}^{T^{N_m}\circ\theta_1}\cal{Z}\circ\Theta_k-\sum_{k=0}^{T^{N_m}-1}\cal{Z}\circ\Theta_k\right)\circ\Theta_n=
\left(\bbm{1}_{\{m\in\cal{C}\}}H\circ\Theta_m\right)\circ\Theta_n=\left(H\bbm{1}_{\{0\in\cal{C}\}}\right)\circ\Theta_{n+m}.
\end{equation*}As an $L^2$-limit,
\begin{equation*}
\lim_{m\to\infty}\ E^{Q^s}\left[\left(H\bbm{1}_{\{0\in\cal{C}\}}\right)\circ\Theta_{n+m}
\,\Bigg|\,\cal{G}_n\right]\overset{\scriptstyle{(\ref{092})}}{=}
\lim_{m\to\infty}\ E^{Q^s}\left[\left(H\bbm{1}_{\{0\in\cal{C}\}}\right)\circ\Theta_{m}
\,\Bigg|\,\cal{G}_0\right]\circ\Theta_n\overset{\scriptstyle{(\ref{126})}}{=}0,
\end{equation*}
thus proving that $M_n$ is a $(\cal{G}_n)$-martingale.\\ 
It now remains to prove (\ref{126}). We consider $B\in L^2(\Gamma^s,
\cal{G}_0,Q^s)$ with $L^2$-norm $\|B\|_2= 1.$ Note that $B$ can be considered as a function of $(w,\lambda_.)$ and 
$(u_m,\omega_m)_{m\leq 0}.$ Then it follows that for fixed $(w,\lambda_.)\in W,$ see (\ref{290}), the random vectors B and
\begin{equation*}
\sum_{k=T^m}^{T^{m+1}-1}\cal{Z}\circ\Theta_k=(X^1_{T^{m+1}}-X^1_{T^m},u_m(T^{m+1}-T^m))-v(T^{m+1}-T^m)\mbox{ for }m\geq 1,
\end{equation*}
are independent under the measure $M^s,$ see (\ref{39}). With these considerations in mind we find that for integer 
$p\geq 1,$
\begin{equation}\label{088}\begin{array}{rcl}
E^{Q^s}\left[\left(H\bbm{1}_{\{0\in\cal{C}\}}\right)\circ\Theta_p \cdot B\right]
&=&\sum_{m\geq 1}E^{Q^s}\left[\left(\sum_{k=T^m}^{T^{m+1}-1}\cal{Z}\circ\Theta_k\right)\cdot B,T^m=p\right]\\[11pt]
&=&\sum_{m\geq 1}E^{P}\left[E^{M^s}\left[\sum_{k=T^m}^{T^{m+1}-1}\cal{Z}\circ\Theta_k\right]
E^{M^s}\left[B\right],T^m=p\right]\\[11pt]
&=& E^{P}\left[\left(E^{M^s}\left[H\right]\bbm{1}_{\{0\in\cal{C}\}}\right)\circ\theta_p\ E^{M^s}\left[B\right]\right].
\end{array}\end{equation}
Then observe that we can find measurable functions $\varphi$ and $\psi$ such that
\begin{equation}\label{131}\begin{array}{l}
E^{M^s}\left[H\right]\bbm{1}_{\{0\in\cal{C}\}}=\varphi\left(T^1,(X^1_t)_{t\geq 0},(\Lambda_n)_
{n\geq 0}\right)\bbm{1}_{\{0\in\cal{C}\}},\\
E^{M^s}\left[B\right]=\psi\left(T^0,(X^1_t)_{t\leq 0},(\Lambda_n)_{n\leq -1}\right),
\end{array}\end{equation}recall the definition of $\Lambda_n$ above (\ref{16}). 
The reason why $E^{M^s}[B]$ depends only on $T^0,\ (X^1_t)_{t\leq 0},\ (\Lambda_n)_{n\leq -1},$
whereas the involved cut times $T^k,k\leq -1,$ are based on the whole trajectory $(X^1_t),\ t\in\bb{R},$ is 
that the information about intersections needed to determine $T^k,k\leq -1,$ can be expressed only by 
 $T^0$ and $(X^1_t), t\leq T^0\leq 0,$ since by definition of $T^0,$ we have that
$(X^1_{(-\infty,k-1]})^R\cap(X^1_{[T^0,\infty)})^R=\emptyset,$ for all $k\leq T^0.$ 
In the sequel we will slightly abuse notation. 
One has to think of the following objects to be defined on an extension of the probability space $(W,\cal{W},P),$ see (\ref{290}) and 
below. Recall that under the measure $P=\bar{P}{\otimes}\bf{\Lambda}^{\varepsilon},$ see (\ref{17}), the 
process $(X^1_t)_{t\in\bb{R}}$ is a two-sided $d_1$-dimensional Brownian motion with $P[X^1_0=0]=1$ which is independent of 
$(\Lambda_n)_{n\in\bb{Z}},$ 
a two-sided sequence of i.i.d. Bernoulli random variables with success parameter $\varepsilon>0,$ see (\ref{92}). 
We are interested in large values of $p$ and set
\begin{equation}
L=\left[\frac{p}{3}\right].
\end{equation}

We introduce a copy $((X^+_t)_{t\in\bb{R}},(\Lambda_j^+)_{j\in\bb{Z}})$ of $((X_t^1)_{t\in\bb{R}},
(\Lambda_j)_{j\in\bb{Z}})$ evolving according to $P$ such that 
$X_t^+=X^1_{t+p}-X^1_p$ for  $t\in [-L,\infty),$ and $\Lambda_j^+=\Lambda_{j+p}$ for $\ j\geq -L,$ and such that 
 $((X^+_t)_{t\in(-\infty,-L)},(\Lambda_j^+)_{j<-L})$ evolves 
independently of $((X^1_{t+p}-X^1_p)_{t\in(-\infty,-L)},(\Lambda_{j+p})_{j<-L}).$ Moreover, we consider another 
copy $((X_t^-)_{t\in\bb{R}},(\Lambda_j^-)_{j\in\bb{Z}})$ of $((X_t^1)_{t\in\bb{R}},(\Lambda_j)_{j\in\bb{Z}})$ 
which is independent of $((X_t^+)_{t\in\bb{R}},(\Lambda_j^+)_{j\in\bb{Z}})$ and evolves according to $P$ 
such that $X^-_t=X^1_t$ for $t\in (-\infty,L],$ and $\Lambda^-_j=\Lambda_j$ for $j\leq L-1,$ and such that 
$((X_t^-)_{t\in(L,\infty)},(\Lambda_j^-)_{j\geq L})$ evolves independently 
of $((X_t^1)_{t\in(L,\infty)},(\Lambda_j)_{j\geq L}).$ Note that \begin{equation}\label{089}
\left((X_t^1)_{t\in\bb{R}},(\Lambda_j)_{j\in\bb{Z}}\right)\overset{\mbox{\scriptsize{law}}}{=}\left((X_t^+)_{t\in\bb{R}},
(\Lambda_j^+)_{j\in\bb{Z}}\right)
\overset{\mbox{\scriptsize{law}}}{=}\left((X_t^-)_{t\in\bb{R}},(\Lambda_j^-)_{j\in\bb{Z}}\right)
\end{equation}and \begin{equation}\label{090}
\left((X_t^+)_{t\in\bb{R}},(\Lambda_j^+)_{j\in\bb{Z}}\right)\mbox{ is independent of }\left((X_t^-)_{t\in\bb{R}},
(\Lambda_j^-)_{j\in\bb{Z}}\right).
\end{equation}The random time $T^-$ is defined like $T^0$ relatively to $((X^-_t)_{t\in\bb{R}},(\Lambda_j^-)_{j\in\bb{Z}})$ and 
$T^+$ is the analogue of $T^1$ attached to $((X^+_t)_{t\in\bb{R}},(\Lambda_j^+)_{j\in\bb{Z}}).$ The random set $\cal{C}^+$ is defined analogously 
to $\cal{C}$ with $((X^+_t)_{t\in\bb{R}},(\Lambda_j^+)_{j\in\bb{Z}})$ in place of $((X_t^1)_{t\in\bb{R}},(\Lambda_j)_{j\in\bb{Z}}),$ 
see (\ref{19}). We then define
\begin{equation}\label{132}\begin{array}{l}
U=E^{M^s}\left[B\right],\ U^-=\psi\left(T^-,(X^-_t)_{t\leq 0},(\Lambda_n^-)_{n\leq -1}\right),\\
V=\left(E^{M^s}\left[H\right]\bbm{1}_{\{0\in\cal{C}\}}\right)\circ\theta_p
=\varphi\left(T^1\circ\theta_p,(X_{t+p}^1-X^1_p)_{t\geq 0},(\Lambda_{n+p})_{n\geq 0}\right)\bbm{1}_{\{p\in\cal{C}\}},

\\ 
V^+=\varphi\left(T^+,(X^+_t)_{t\geq 0},(\Lambda_n^+)_{n\geq 0}\right)\bbm{1}_{\{0\in\cal{C}^+\}}.
\end{array}\end{equation}
By construction, see in particular (\ref{089}) and (\ref{090}), we have that 
$U\overset{\mbox{\tiny{law}}}{=}U^-$ and due to the invariance of $P$ under the shift
$\theta_p,$ also $V\overset{\mbox{\tiny{law}}}{=}V^+,$ but 
$U^-$ and $V^+$ are now independent. For $p\geq 1,$
\begin{eqnarray*}
E^{Q^s}\left[\left(H\bbm{1}_{\{0\in\cal{C}\}}\right)\circ\Theta_p\cdot B\right]&\overset{\scriptstyle{(\ref{088})}}{=}&E^{P}\left[VU\right]
\\&=&E^{P}\left[V^+U^-\right]+E^{P}\left[V^+(U-U^-)\right]+E^{P}\left[(V-V^+)U\right].
\end{eqnarray*}
Note that the first term in the last line vanishes because of the independence mentioned above and the fact that
\begin{equation*}
E^{P}\left[V^+\right]=E^{P}\left[V\right]\overset{\mbox{\scriptsize{(\ref{088})}}}{=}
E^{Q^s}\left[H\bbm{1}_{\{0\in\cal{C}\}}\right]\overset{\mbox{\scriptsize{(\ref{24})}}}{=}
E^{\hat{Q}^s}\left[H\right]{E^{\hat{P}}}\left[T^1\right]^{-1}
\overset{{(\ref{234})}}{=}\ E^{Q^s}\left[\cal{Z}\right]\overset{\scriptstyle{(\ref{121})}}{=}0.
\end{equation*}
Therefore, after recalling that 
$\|U\|_2\leq\|B\|_2= 1,$ we find with H\"older's inequality:
\begin{equation}\label{133}
E^{Q^s}\left[\left(H\bbm{1}_{\{0\in\cal{C}\}}\right)\circ\Theta_p\cdot B\right]\leq
\|V^+\|_4\|U-U^-\|_{4/3}+\|V-V^+\|_2.
\end{equation}
Due to stationarity of $\theta_1$ under $P$ and Jensen's inequality we easily obtain that
\begin{equation*}
\|V^+\|_4=\|V\|_4\leq E^{Q^s}\left[|H|^4\bbm{1}_{\{0\in\cal{C}\}}\right]^{1/4}=
 E^{\hat{Q}^s}\left[|H|^4\right]^{1/4}P\left[0\in\cal{C}\right]^{1/4}\leq
E^{\hat{Q}^s}\left[|H|^4\right]^{1/4}.
\end{equation*}From the definition of $H,$ see (\ref{127}), and Remark (\ref{540}) it then follows that\begin{equation}\label{134}
\|V^+\|_4=\|V\|_4\leq E^{\hat{P}\times K_{0}}\left[|\chi_{T^1}|^4\right]^{1/4}+vE^{\hat{P}}
\left[(T^1)^4\right]^{1/4}\overset{{\scriptstyle{(\ref{600}), (\ref{602})}}}{<}\infty.
\end{equation}
In view of the definitions (\ref{132}) we see that
\begin{equation}\label{050}
\|V-V^+\|_2\leq \left\|\left(|V|+|V^+|\right)\left(\bbm{1}_{\{T^+\ne T^1\circ\theta_p\}}+
|\bbm{1}_{\{p\in\cal{C}\}}-\bbm{1}_{\{0\in\cal{C}^+\}}|\right)\right\|_2.
\end{equation}
Since by stationarity of $\theta_1$ under $P$ and the identity in law (\ref{089}), 
$\ P\left[\{p\in\cal{C}\}\smallsetminus\{0\in\cal{C}^+\}\right]$ is equal to
$P\left[\{0\in\cal{C}^+\}\smallsetminus\{p\in\cal{C}\}\right],$ an application of Cauchy-Schwarz' inequality to the 
right-hand side of (\ref{050}) shows that
\begin{equation}\label{051}
\|V-V^+\|_2\leq 2\|V\|_4\left(P\left[T^+\ne T^1\circ\theta_p\right]^{1/4}+2P\left[\{p\in\cal{C}\}\smallsetminus
\{0\in\cal{C}^+\}\right]^{1/4}\right).
\end{equation}
Since $(X^+_t,\Lambda^+_n)$ and $(X^1_t,\Lambda_n)\circ\theta_p$ coincides for $t\in [-L,\infty),\ n\geq -L,$ 
with (\ref{19}) we see that for large $p,$ the events $\{T^+\ne T^1\circ\theta_p\}$ and $\{p\in\cal{C}\}\smallsetminus\{0\in\cal{C}^+\}$ are both included in 
\begin{equation*}
\left\{\left(X^+_{(-\infty,-L]}\right)^{R}\cap\left(X^+_{[0,\infty)}\right)^{R}\ne\emptyset\right\}\cup
\left\{\left(\left(X^1_.\circ\theta_p\right)_{(-\infty,-L]}\right)^{R}\cap\left(\left(X^1_.\circ\theta_p\right)_
{[0,\infty)}\right)^{R}\ne\emptyset\right\}, 
\end{equation*}and so, together with (\ref{051}) we find using stationarity once again that
\begin{equation}\label{135}
\|V-V^+\|_2\leq c\|V\|_4P\left[\left(X^1_{(-\infty,0]}\right)^{R}\cap \left(X^1_{[L,\infty)}\right)^{R}
\ne\emptyset\right]^{1/4}.
\end{equation}
By analogous arguments as above we also find that,
\begin{equation}\label{136}\begin{array}{rcl}
&&\|U-U^-\|_{4/3} \\[5pt]
&&\leq  \left\|\left(|U|+|U^-|\right)\bbm{1}_{\{T^0\ne T^-\}}\right\|_{4/3}\\[5pt]
&&\leq \Big\|\left(|U|+|U^-|\right)\Big(\bbm{1}
{\scriptscriptstyle{\left\{\left(X^-_{(-\infty,0]}\right)^{R}\cap\left(X^-_{[L,\infty)}\right)^{R}\ne\emptyset\right\}}}
 +\bbm{1}
{\scriptscriptstyle\left\{\left(X^1_{(-\infty,0]}\right)^{R}\cap\left(X^1_{[L,\infty)}\right)^{R}\ne\emptyset\right\}}\Big)
\Big\|_{4/3}\\
&& \leq  cP\left[\left(X^1_{(-\infty,0]}\right)^{R}\cap \left(X^1_{[L,\infty)}\right)^{R}
\ne\emptyset\right]^{1/4},
\end{array}\end{equation}
where we used H\"older's inequality and $\|U^-\|_2=\|U\|_2\leq\|B\|_2=1$ in the last inequality. Collecting (\ref{133}), (\ref{134}), 
(\ref{135}) and  (\ref{136}), we finally find\begin{eqnarray*}
\left\|E^{Q^s}\left[\left(H\bbm{1}_{\{0\in\cal{C}\}}\right)\circ\Theta_p\ |\ \cal{G}_0\right]\right\|_2 & \leq &
c\|V\|_4P\left[\left(X^1_{(-\infty,0]}\right)^{R}\cap \left(X^1_{[L,\infty)}\right)^{R}
\ne\emptyset\right]^{1/4}\\
&\overset{\scriptstyle{(\ref{33})}}{\leq}&c\|V\|_4\ p^{-\frac{d_1-4}{8}}.
\end{eqnarray*}This quantity is summable in $p,$ since $d_1\geq 13.$ This finishes the proof of (\ref{126}) and thus 
of Theorem \ref{120}.
\begin{flushright}$\Box$\end{flushright}
\pagebreak

{\rem{
In the next section we will strengthen Theorem \ref{103} and \ref{120} to central limit theorems under the quenched measure. In the literature only 
few results on quenched invariance principles for diffusions in random environment are available. One result is due to 
Sznitman-Zeitouni \cite{SZNZEIT} who consider small perturbations of Brownian motion, and a second situation in which a 
quenched central limit theorem holds is discussed in Osada \cite{O}. The latter result is attained with the technique 
of the {\textit{environment viewed from the particle.}}
}}

\section{Central limit theorem under the quenched measure}\label{d}

We are going to show how one can improve the results of Section \ref{c} to central limit theorems under the 
quenched measure $P_{0,\omega}$. We use an idea of Bolthausen and Sznitman, see Lemma 4 in \cite{BS}, to turn the {\textit{annealed}} invariance principle into a {\textit{quenched}} invariance principle,  
by bounding certain variances through the control of intersections of two independent paths. For this 
purpose we do not require an explicit invariant measure for the process of the environment viewed from the particle or 
the control of moments of certain regeneration times, see for instance in \cite{berzei}, \cite{RS2}, \cite{RS4}, in 
the discrete setting. We recall the definition of $v$ in (\ref{55}).\\

{\thm{\label{200}Assume $d_1\geq 7$ and (\ref{100}), or $d_1\geq 13.$ Then for $\bb{P}$-a.e.$\,\omega,$ under the measure $P_{0,\omega},$ 
the $C(\bb{R}_+,\bb{R}^d)$-valued random variables\begin{equation*}
B^r_.\df\frac{1}{\sqrt{r}}\left(X_{r\cdot}-vr\cdot\right),\qquad r>0,
\end{equation*} converge weakly to a Brownian motion $B_.$ with covariance matrix $A$ given in Theorem \ref{103} and 
\ref{120} respectively, as $r\to\infty.$}}\\

{\textbf{Proof:}} By similar arguments as at the beginning of the proof of Theorem \ref{103}, see (\ref{301}) and 
(\ref{321}), and the identity in law in (1) of Theorem \ref{47} we can see that it suffices to show that
\begin{equation}\label{161}\begin{array}{l}\mbox{for }\bb{P}\mbox{-a.e. }\omega,
\mbox{ under the measure }P\times K_{0,\omega},\mbox{ the }C(\bb{R}_+,\bb{R}^d)\mbox{-valued\hspace{2pt} random}\\
\mbox{variables }\beta^n_.\df\frac{1}{\sqrt{n}}\left\{\chi_{[n\cdot]}^{}+\left(n\cdot-[n\cdot]\right)\left(
\chi_{[n\cdot]+1}^{}-\chi_{[n\cdot]}^{}\right)-vn\cdot\right\},\mbox{\,for\,integers}\\
n\geq 1,\mbox{ converge weakly to }B_.,\mbox{ as }n\to\infty.
\end{array}\end{equation}
From the proofs of Theorem \ref{103} and \ref{120}, see in particular (\ref{301}), (\ref{321}), (\ref{060}) and (\ref{122}),
we know that \begin{equation}\label{061}
\beta_.^n\longrightarrow B_.\mbox{ in law under }P\times K_0,\mbox{ as }n\to\infty.
\end{equation}
From the proof of Lemma 4.1 in \cite{BS} we see that (\ref{161}) follows from (\ref{061}) and  a variance calculation.
Let us introduce for $T>0$ the space of continuous, $\bb{R}^d$-valued functions on 
$[0,T]$ denoted with $C([0,T],\bb{R}^d),$ 
which we equip with the distance\begin{equation}\label{304}
d_T(g,g')\df\sup_{t\leq T}|g(t)-g'(t)|\wedge 1.
\end{equation}The proof of Lemma 4.1 in \cite{BS} shows us that (\ref{161}) follows once we prove that for all 
$T>0,\ \xi\in (1,2]$ and all bounded 
Lipschitz functions $F$ on $C([0,T],\bb{R}^d),$\begin{equation}\label{160}
\sum_m \mathrm{Var}_{\bb{P}}\left(E^{P\times K_{0,\omega}}\left[F\left(\beta_.^{[\xi^m]}\right)\right]\right)<\infty
\end{equation} (with a slight abuse of notation). For this purpose we need some further notation. Given an environment 
$\omega,$ we consider two independent copies $((\chi_t^{})_{t\geq 0}^{},(\Lambda_n)_{n\geq 0}^{})$ and 
$((\tilde{\chi}_t^{})_{t\geq 0}^{},(\tilde{\Lambda}_n)_{n\geq 0}^{})$ evolving according to $P\times K_{0,\omega}.$ 
The corresponding first $d_1$ components of 
$\chi_.^{}$ and $\tilde{\chi}_.^{}$ denoted with 
 $X^1_.$ and $\tilde{X}_.^1$ are then two independent $d_1$-dimensional Brownian motions. We also 
introduce the corresponding polygonal interpolations $\beta_.^n$ and $\tilde{\beta}_.^n$ defined as in (\ref{161}). With $\cal{D}$ we denote the 
set of one-sided cut times attached to $((X^1_t)_{t\in\bb{R}},(\Lambda_j)_{j\in\bb{Z}})$ defined via \begin{equation}\label{091}
\cal{D}=\left\{k\geq 1\ \Big|\ \left(X^1_{[0,k-1]}\right)^{R}{{\cap}}
\left(X^1_{[k,\infty)}\right)^{R}=\emptyset\mbox{ and }\Lambda_{k-1}=1\right\}.
\end{equation}$\tilde{\cal{D}}$ is defined analogously and attached to $((\tilde{X}_t^1)_{t\in\bb{R}},
(\tilde{\Lambda}_j)_{j\in\bb{Z}}).$ We then pick
\begin{equation*}\xi\in (1,2],\quad 0<\mu<\nu<\frac{1}{2},\end{equation*}and for $m\geq 1$ we define $n=[\xi^m],$ 
as well as
\begin{equation*}\rho_m\df\inf\left\{\cal{D}\,{{\cap}}\,[n^{\mu},\infty)\right\}<\infty,\quad P\mbox{-a.s.
 (see (\ref{22}))},\end{equation*} and $\tilde{\rho}_m$ as the corresponding variable attached to $((\tilde{X}_t^1)_{t\in\bb{R}},
(\tilde{\Lambda}_j)_{j\in\bb{Z}}).$
 In order to take advantage of decoupling effects we will consider the event\begin{equation*}
 \cal{A}_m=\left\{\rho_m\vee\tilde{\rho}_m\leq n^{\nu},\ \left(X^1_{[0,\infty)}\right)^{R}{{\cap}}
\left(\tilde{X}^1_{[n^{\mu},\infty)}\right)^{R}=\emptyset,\ \left(X^1_{[n^{\mu},\infty)}\right)^{R}{{\cap}}
\left(\tilde{X}^1_{[0,\infty)}\right)^{R}=\emptyset\right\}.
 \end{equation*}
 We are now ready to prove (\ref{160}). Without loss of generality, we assume the Lipschitz constant and the absolute 
value of $F$ to be bounded 
 by 1. For the remainder of the proof we write $E$ and $E_{\omega}$ for the expectation under the measure $P\times K_0$ and 
 $P\times K_{0,\omega}$ respectively. For $m\geq 1,$ we have
 \begin{eqnarray*}
 {\mathrm{Var}}_{\bb{P}}\left(E_{\omega}\left[F\left(\beta_.^n\right)\right]\right)&=&
 \bb{E}\left[E_{\omega}{\otimes}E_{\omega}\left[F\left(\beta_.^n\right)F(\tilde{\beta}_.^n)
\right]\right]-E{{\otimes}}E\left[F\left(\beta_.^n\right)F(\tilde{\beta}_.^n)\right]\\
&=& \bb{E}\left[E_{\omega}{{\otimes}}E_{\omega}\left[F\left(\beta_.^n\right)F(\tilde{\beta}_.^n),
\cal{A}_m\right]\right]-E{{\otimes}}E\left[F\left(\beta_.^n\right)F(\tilde{\beta}_.^n),
\cal{A}_m\right]+d_m
\end{eqnarray*} with \begin{equation}\label{142}|d_m|\leq 2P{{\otimes}}P\left[\cal{A}_m^c\right]. 
\end{equation} Using that $F$ is bounded and Lipschitz and $d_T(\cdot,\cdot)\leq 1$ we obtain that the difference of the first 
two terms in the last line above (\ref{142}) is equal to\begin{equation}\label{143}\begin{array}{l}
\bb{E}\left[E_{\omega}{{\otimes}}E_{\omega}\left[F\left(\beta_{\cdot+\frac{\rho_m}{n}}^n-\beta^n_{\frac{\rho_m}{n}}\right)
F\left(\tilde{\beta}_{\cdot+\frac{\tilde{\rho}_m}{n}}^n-\tilde{\beta}^n_{\frac{\tilde{\rho}_m}{n}}\right),\cal{A}_m\right]\right]\\[11pt]
-E{{\otimes}}E\left[F\left(\beta_{\cdot+\frac{\rho_m}{n}}^n-\beta^n_{\frac{\rho_m}{n}}\right)
F\left(\tilde{\beta}_{\cdot+\frac{\tilde{\rho}_m}{n}}^n-\tilde{\beta}^n_{\frac{\tilde{\rho}_m}{n}}\right),\cal{A}_m\right]
+\Delta^3_m\\[11pt]
 =:\Delta_m^1-\Delta^2_m+\Delta^3_m,\end{array}\end{equation}
with \begin{equation*}
\Delta^3_m\leq 6E{{\otimes}}E\left[d_T\left(\beta_{\cdot+\frac{\rho_m}{n}}^n-\beta^n_{\frac{\rho_m}{n}}\,,\,
\beta_.^n\right),\cal{A}_m\right].
\end{equation*}We first want to show that\begin{equation}\label{805}\Delta_m^1=\Delta_m^2.
\end{equation}For each $\omega\in\Omega$ and fixed samples $(w,\lambda_.)$ and $(\tilde{w},\tilde{\lambda}_.)$ of 
$(X^1_.,\Lambda_.)$ and $(\tilde{X}^1_.,\tilde{\Lambda}_.)$ respectively, 
under $K_{0,\omega}(w,\lambda_.){{\otimes}}
K_{0,\omega}(\tilde{w}_.,\tilde{\lambda}_.),$ the processes $\beta^n_.$ and $\tilde{\beta}_.^n$ 
are independent and hence with the help of Fubini's Theorem we can write
\begin{equation*}
\Delta_m^1=E^{P}{{\otimes}}E^{P}\left[\bb{E}\left[E^{K_{0,\omega}}\left[
F\left(\beta_{\cdot+\frac{\rho_m}{n}}^n-\beta^n_{\frac{\rho_m}{n}}\right)\right]
E^{K_{0,\omega}}\left[F\left(\tilde{\beta}_{\cdot+\frac{\tilde{\rho}_m}{n}}^n-\tilde{\beta}^n_{\frac{\tilde{\rho}_m}{n}}\right)\right]\right],
\cal{A}_m\right].
\end{equation*}By similar arguments to those leading to (\ref{701}) we obtain that for each $(w,\lambda_.)\in W,$
\begin{equation}\label{071}\begin{array}{l}
E^{K_{0,\omega}}\left[F\left(\beta_{\cdot+\frac{\rho_m}{n}}^n-\beta^n_{\frac{\rho_m}{n}}\right)\right]\\[11pt]
={\displaystyle{\int_{\bb{R}^{d_2}}\frac{dy}{vol(d_2)}}}E^{K_{0,\omega}}\left[\bbm{1}{\{y\in B^{d_2}_1(X^2_{\rho_m-1})\}}\right]
E^{K_{0,\bar{\omega}}\circ\theta_{\rho_m}}\Big[F\left(\beta^n_.-\beta^n_0\right)\circ\theta_{\rho_m}\Big],
\end{array}\end{equation}with $\bar{\omega}=\tau_{(w(\rho_m),y)}(\omega).$ From (4) of Theorem \ref{47} it follows that the first 
expectation in the second line of (\ref{071}) is measurable with respect to $\cal{H}_{(w([0,\rho_m-1]))\times\bb{R}^{d_2}},$ see (\ref{4}),
whereas the second expectation is $\cal{H}_{(w([\rho_m,\infty)))\times\bb{R}^{d_2}}$-measurable. With these considerations in mind
we find with Fubini's Theorem and finite range dependence, see (\ref{5}), that on $\cal{A}_m,$
\begin{eqnarray}
\label{070} && \bb{E}\left[E^{K_{0,\omega}}\left[F\left(\beta_{\cdot+\frac{\rho_m}{n}}^n-\beta^n_{\frac{\rho_m}{n}}\right)\right]
E^{K_{0,\omega}}\left[F\left(\tilde{\beta}_{\cdot+\frac{\tilde{\rho}_m}{n}}^n-\tilde{\beta}^n_{\frac{\tilde{\rho}_m}{n}}\right)\right]\right]\\[5pt]
\nonumber && =\int_{\bb{R}^{d_2}}\int_{\bb{R}^{d_2}}\frac{dy_1dy_2}{vol(d_2)^2}
 \bb{E}\left[E^{K_{0,\omega}}\left[\bbm{1}{\{y_1\in B^{d_2}_1(X^2_{\rho_m-1})\}}\right]
E^{K_{0,\omega}}\left[\bbm{1}{\{y_2\in B^{d_2}_1(\tilde{X}^2_{\tilde{\rho}_m-1})\}}\right]\right]\\[5pt]
\nonumber &&\hspace{2cm}\times \bb{E}\left[E^{K_{0,\bar{\omega}}\circ\theta_{\rho_m}}\Big[F\Big(\beta^n_.-\beta^n_0\Big)\circ\theta_{\rho_m}\Big]\right]
 \bb{E}\left[E^{K_{0,\tilde{\omega}}\circ\theta_{\tilde{\rho}_m}}\Big[F\left(\tilde{\beta}^n_.-\tilde{\beta}^n_0\right)\circ\theta_{\tilde{\rho}_m}\Big]\right],
\end{eqnarray}with $\bar{\omega}=\tau_{(w(\rho_m),y_1)}(\omega),\ \tilde{\omega}=\tau_{(w(\tilde{\rho}_m),y_2)}(\omega).$ Because of the 
stationarity of the environment the last two $\bb{P}$-expectations above are in fact independent of $y_1$ respectively $y_2$ so that an 
application of Fubini's Theorem shows us that (\ref{070}) equals\begin{equation*}
E^{K_{0}\circ\theta_{\rho_m}}\Big[F\Big(\beta^n_.-\beta^n_0\Big)\circ\theta_{\rho_m}\Big]
 E^{K_{0}\circ\theta_{\tilde{\rho}_m}}\Big[F\left(\tilde{\beta}^n_.-\tilde{\beta}^n_0\right)\circ\theta_{\tilde{\rho}_m}\Big].
\end{equation*} Analogously we also find that\begin{equation*}
E^{K_{0}}\Big[F\Big(\beta^n_{\cdot+\frac{\rho_m}{n}}-\beta^n_{\frac{\rho_m}{n}}\Big)\Big]=
E^{K_{0}\circ\theta_{\rho_m}}\Big[F\Big(\beta^n_.-\beta^n_0\Big)\circ\theta_{\rho_m}\Big],
\end{equation*}and the same holds true if we replace $\beta^n_.,\rho_m$ by $\tilde{\beta}^n_.,\tilde{\rho}_m.$ 
This concludes the proof of (\ref{805}). We now come to the control of $\Delta_m^3.$ Noting that 
on $\cal{A}_m,\ E_{\omega}$-a.s.\begin{equation*}
d_T\left(\beta_{\cdot+\frac{\rho_m}{n}}^n-\beta^n_{\frac{\rho_m}{n}}\,,\,\beta_.^n\right)\leq
\sup_{\scriptscriptstyle{\substack{0\leq s< t\leq Tn+1\\ \,|t-s|\leq n^{\nu}}}}\frac{1}{\sqrt{n}}\big|
\chi_t^{}-\chi_s^{}\big|+\sup_{\scriptscriptstyle{t\leq n^{\nu}}}\frac{1}{\sqrt{n}}\big|\chi_t^{}\big|,
\end{equation*}
we find by using (\ref{150}) and the fact that \begin{equation*}
E_{\omega}\Bigg[\sup_{\scriptscriptstyle{\substack{0\leq s< t\leq Tn+1\\ \,|t-s|\leq n^{\nu}}}}
\big|W_t-W_s\big|\Bigg]\leq 
c(T)n^{1/4+\nu/4}E_{\omega}\Bigg[\sup_{0\leq s\leq t\leq 1}\frac{|W_t-W_s|}{|t-s|^{1/4}}\Bigg]
\leq c(T)n^{1/4+\nu/4},
\end{equation*}where the last inequality follows from an application of Fernique's Theorem, see \cite{DS} on page 14, which ensures the 
existence of the exponential moment of $\eta\sup_{0\leq s\leq t\leq 1}|W_t-W_s|/|t-s|^{1/4}$ for a certain constant $\eta>0,$
that\begin{equation*}
\Delta^3_m\leq \frac{c(T)}{\sqrt{n}}(n^{\nu}+n^{1/4+\nu/2})\leq c(T)n^{\nu/2-1/4},
\end{equation*} and hence $\sum_m \Delta^3_m<\infty,$ recall $n=[\xi^m].$ It remains to show that\begin{equation}\label{144}
\sum_m P{{\otimes}}P\left[\cal{A}_m^c\right]<\infty.
\end{equation}
Indeed, we find that\begin{equation}\label{145}
P{{\otimes}}P\left[\left(X^1_{[0,\infty)}\right)^{R}{{\cap}}
\left(\tilde{X}^1_{[n^{\mu},\infty)}\right)^{R}\ne\emptyset\right]\overset{\scriptstyle{(\ref{33})}}{\leq} cn^{-\mu\frac{d_1-4}{2}},
\end{equation}
and moreover, since the random set $\cal{C}\,{{\cap}}\,\bb{N}$ is contained in $\cal{D},$ see (\ref{19}) and 
(\ref{091}), we have that $P$-a.s., $\rho_m-n^{\mu}\leq T^1\circ\theta_{[n^{\mu}]}$ and hence from stationarity of 
$\theta_1$ under $P$ it follows that for large $m,$
\begin{equation}\label{146}
P\left[\rho_m>n^{\nu}\right]\leq P[T^1>n^{\nu}-n^{\mu}]\overset{\scriptstyle{(\ref{26})}}{\leq}
c(\varepsilon)(\log n^{\nu})^{1+\frac{d_1-4}{2}}n^{-\nu\frac{d_1-4}{2}}\leq e^{-c(\varepsilon)m}.
\end{equation} Combining (\ref{145}) and (\ref{146}) we deduce (\ref{144}).
\begin{flushright}$\Box$\end{flushright}

\section\appendixname\label{e}
\subsection{Two central limit theorems for martingales}
For $T>0$ and integer $d\geq 1$ we denote with $D([0,T],\bb{R}^d)$ the space of all 
$\bb{R}^d$-valued c$\grave{\mathrm{a}}$dl$\grave{\mathrm{a}}$g functions on $[0,T]$ which we endow with the Skorohod metric. The 
space of all continuous functions on $[0,T]$ with values in $\bb{R}^d$ is denoted with 
$C([0,T],\bb{R}^d)$ and is equipped with the supnorm. See Chapter 2 and 3 of \cite{B} 
for an extensive discussion of the above mentioned function spaces.

{\thm{\label{080}
$X_n,\ \cal{F}_n\df\sigma(X_k,\ k\leq n),\ n\geq 1,$ is an $\bb{R}^{d}$-valued sequence of square integrable martingale 
differences on a probability space $(\Omega,\cal{F},P),$ i.e. $E[|X_n|^2]<\infty$ and $E[X_n|\cal{F}_{n-1}]=0.$ Let $\Gamma$ 
be a symmetric, non-negative definite $d\times d$-matrix and $S_n(\cdot)\df\sum_{k=1}^{[n\cdot]}X_k.$ Assume that
\begin{equation}\label{081}
\lim_{n\to\infty}\frac{1}{n}\sum_{k=1}^{[ns]}E\left[X_k X_k^t\,|\,\cal{F}_{k-1}\right]=s\Gamma\quad\mbox{in probability,}
\end{equation} for each $s\in\bb{R}_+,$ and\begin{equation}\label{082}
\lim_{n\to\infty}\frac{1}{n}\sum_{k=1}^n E\left[|X_k|^2 \bbm{1}\left\{|X_k|\geq \varepsilon\sqrt{n}\right\}\,|\,\cal{F}_{k-1}\right]=0\quad
\mbox{in probability,}
\end{equation}for each $\varepsilon>0.$ Then the $C(\bb{R}_+,\bb{R}^d)$-valued random variables\begin{equation}\label{083}
\frac{1}{\sqrt{n}}\bar{S}_n(\cdot)\df\frac{1}{\sqrt{n}}\left\{S_n(\cdot)+\left(n\cdot-[n\cdot]\right)\left(
X_{[n\cdot]+1}-X_{[n\cdot]}\right)\right\},\quad n\geq 1,
\end{equation}converge weakly to a $d$-dimensional Brownian motion $B_.$ with covariance matrix $\Gamma,$ as $n\to\infty.$
}}\\

{\textbf{Proof:}} An application of an invariance principle for a vector-valued, square-integrable 
martingale difference array which is proved in \cite{RS}, see Theorem 3, shows us that\begin{equation*}
\frac{1}{\sqrt{n}}S_n(\cdot)\mbox{ converges weakly to }B_.\mbox{ on the Skorohod space }D([0,1],\bb{R}^d),\mbox{ as }n\to\infty.
\end{equation*}From (18) in the proof of Theorem 3 in \cite{RS} we know that in view of Remark \ref{302} also the polygonal interpolation 
$n^{-1/2}\bar{S}_{n\cdot}$ converges weakly on $D([0,1],\bb{R}^{d})$ to the same limit. Since the Skorohod topology 
relativized to $C([0,1],\bb{R}^{d})$ coincides with 
the uniform topology there we have in fact weak convergence on $C([0,1],\bb{R}^{d})$ for the process 
$n^{-1/2}\bar{S}_{n\cdot}.$ Moreover, the identity\begin{equation*}
\frac{1}{\sqrt{n}}\bar{S}_{n\cdot}=\sqrt{M}\frac{1}{\sqrt{nM}}\bar{S}_{nM\frac{\cdot}{M}},
\qquad M\geq 1,
\end{equation*}shows that the process $n^{-1/2}\bar{S}_{n\cdot}$ indeed converges 
weakly on $C([0,M],\bb{R}^{d}).$ Then weak convergence on 
each $C([0,M],\bb{R}^{d})$ implies weak convergence on $C(\bb{R}_+,\bb{R}^{d}),$ see \cite{W}.\begin{flushright}$\Box$\end{flushright}

{\thm{\label{140}$X_n,\ \cal{F}_n\df\sigma (X_m, m\leq n),\ n\in\bb{Z},$ is an $\bb{R}^d$-valued, ergodic 
stationary sequence of square integrable martingale differences on a
probability space $(\Omega,\cal{F},P)$, 
i.e. $E[|X_n|^2]=E[|X_1|^2]<\infty$ and $E[X_n|\cal{F}_{n-1}]=0.$ Let 
$\Gamma\df E[X_1X_1^t]$ and $S_n(\cdot)\df\sum_{k=1}^{[n\cdot]}X_k.$ Then, the $C(\bb{R}_+,\bb{R}^d)$-valued random variables  $$
\frac{1}{\sqrt{n}}\bar{S}_n(\cdot)\df \frac{1}{\sqrt{n}}\left\{S_n(\cdot)+(n\cdot-[n\cdot])\left(X_{[n\cdot]+1}-X_{[n\cdot]}\right)\right\},
\quad n\geq 1,$$ 
converge weakly to a 
$d$-dimensional Brownian motion $B_.$ with covariance matrix $\Gamma,$ as $n\to\infty.$}}\\

\textbf{Proof:} Note that for an $\bb{R}^m$-valued function $f$ with $m\geq 1$ 
such that $f(X_k)\in L^1(P), k\in\bb{Z},$ the conditional expectation $E[f(X_k)\,|\,\cal{F}_{k-1}]$ 
can be written as $\varphi(X_{k-1},X_{k-2},\ldots)$ for a measurable function $\varphi,$ which does not depend on $k.$ 
By the ergodic theorem we thus find that\begin{equation*}
P\mbox{-a.s.},\ \lim_{n\to\infty}\frac{1}{{n}}\sum_{k=1}^{[ns]}E\left[X_k X_k^t\,|\,\cal{F}_{k-1}\right]=
sE\left[X_1X_1^t\right]=s\Gamma.
\end{equation*}
Let $\delta>0$ be small and choose $N=N(\delta)$ such that\begin{equation*}
E\left[|X_1|^2,|X_1|\geq\varepsilon\sqrt{N}\right]<\delta.
\end{equation*}Then,\begin{eqnarray*}
\lim_{n\to\infty} \frac{1}{{n}}\sum_{k=1}^n E\left[|X_k|^2\bbm{1}_{\{|X_k|\geq\varepsilon\sqrt{n}\}}\,\big|\,\cal{F}_{k-1}\right]
&\leq & \lim_{n\to\infty}\frac{1}{n}\sum_{k=1}^n E\left[|X_{k}|^2\bbm{1}_{\{|X_{k}|\geq\varepsilon\sqrt{N}\}}\,\big|
\,\cal{F}_{k-1}\right]\\
& = & E\left[|X_1|^2,|X_1|\geq\varepsilon\sqrt{N}\right]<\delta,\quad P\mbox{-a.s.},
\end{eqnarray*}
where we used the ergodic theorem in the last equality. An application of Theorem \ref{080} concludes the proof.
\begin{flushright}$\Box$\end{flushright}

\end{document}